\documentclass[12pt]{article}

\usepackage{amsmath,amssymb,amsthm}
\usepackage{verbatim}
\usepackage{color}

\def\indicator{{\mathchoice {1\mskip-4mu\mathrm l}%
{1\mskip-4mu\mathrm l}{1\mskip-4.5mu\mathrm l}%
{1\mskip-5mu\mathrm l}}}

\newtheorem{theorem}{Theorem}[section]
\newtheorem{corollary}[theorem]{Corollary}
\newtheorem{lemma}[theorem]{Lemma}
\newtheorem{proposition}[theorem]{Proposition}

\newtheorem{remark}[theorem]{Remark}
\numberwithin{equation}{section}

\newtheorem{assumption}{Assumption}

\newtheorem{assumptionpr}{Assumption}

\def\be{\begin{equation}}
\def\ee{\end{equation}}

\def\q{\quad}
\def\qq{{\qquad}}

\def\wh{\widehat}

\def\eps{\varepsilon} 

\def\phi{\varphi}

\def\om{{\omega}}
\def\ol{\overline}

\def\proof{{\medskip\noindent {\bf Proof. }}}

\def\qed{{\hfill $\square$ \bigskip}}

\def\dist{{\mathop {{\rm dist\, }}}}

\def\cov{{\mathop {{\rm cov }}}}

  \def\sC {{\cal C}}
 \def\sE {{\cal E}} 
  
  \def\sL {{\cal L}}

  \def\sX {{\cal X}}

 \def\bE {{\mathbb E}}

 \def\bN {{\mathbb N}} 
\def\bP {{\mathbb P}}  \def\bR {{\mathbb R}}

 \def\bZ {{\mathbb Z}}

\def\ignore#1{}  

\def\ms{\medskip}

\def\sm{\smallskip\noindent}

\def\obE{\ol \bE}

\begin{document}

\title{\bf Invariance Principle for the Random Conductance Model with dynamic bounded Conductances}

\author{Sebastian Andres \footnote{Research partially supported by NSERC (Canada)}
\date{}
}

\maketitle

\begin{abstract}
We study a continuous time random walk $X$ in an environment of dynamic
random conductances in $\mathbb{Z}^d$.
We assume that the conductances are stationary ergodic, uniformly bounded and bounded away from zero and polynomially mixing in space and time.
We prove a quenched invariance principle for $X$, 
and obtain Green's functions bounds and a local limit theorem. We also discuss a connection to stochastic interface models.

\end{abstract}

\section{Introduction} \label{sec:intro}

We consider the Euclidean lattice $\bZ^d$ equipped with the set $E_d$ of non oriented nearest neighbour bonds:
$E_d=\{e=\{x,y\}: x,y\in\bZ^d, |x-y|=1\}$. We will also write $x\sim y$ when $\{x,y\} \in E_d$. Denote by $\hat\Omega=[0,\infty)^{E_d}$ and by $\Omega$ 
 the set of all measurable functions from $\bR$ to $\hat \Omega$. We equip $\Omega$ with a $\sigma$-algebra $\mathcal{F}$ and a probability measure $\bP$ so that $(\Omega, \mathcal{F},\mathbb{P})$ becomes a probability space.  The random environment is given by the coordinate maps $\mu_{e}^\om(t)=\om_e(t)$, $t\in \mathbb{R}, e\in E_d$. We will refer to $\mu_e(t)$ as the \emph{conductance} of the edge $e$ at time $t$. Further,  write $\mu^\om_{xy}(t)=\mu_{\{x,y\}}(t)=\mu_{yx}(t)$, 
and $\mu_{xy}(t)=0$ if $\{x, y\} \not\in E_d$, and set 
\begin{align}   \label{jumpP}
  \mu_x(t) = \sum_{y \in\bZ^d}  \mu_{xy}(t)=\sum_{y \sim x}  \mu_{xy}(t). 
\end{align}

We denote by $D(\bR,\bZ^d)$ the space of  $\bZ^d$-valued c\`{a}dl\`{a}g functions on $\bR$. For a given $\om\in\Omega$ and for $s\in \bR$ and $x\in\bZ^d$, let  $P_{s,x}^\omega$ be the probability measure on $D(\bR,\bZ^d)$, under which the coordinate process $(X_t)_{t\in\bR}$ is the continuous-time Markov chain on $\bZ^d$ starting in $x$ at time $t=s$ with time-dependent generator given by:
\begin{align}\label{e-LV}
 \sL^\om_t f(x) = \sum_{y\sim x} \mu^\om_{xy}(t) (f(y)-f(x)).
\end{align}
That is, $X$ is the time-inhomogeneous random walk, whose time-dependent jump rates are given by the conductances.
Note that the counting measure, independent of $t$, is an invariant measure for $X$. 
Further, we denote by $p^\om(s,x;t,y)$, $x,y\in \bZ^d$, $s\leq t$, the transition densities of the time-inhomogeneous random walk $X$.
This model of a random walk in a random
environment is  known in the literature -- at least in the case of time-independent conductances -- as  the {\em Random Conductance Model} or RCM. Note that the total jump rate out of any site $x$ is not normalized, in particular the sojourn time at site $x$ depends on $x$. Therefore, the random walk $X$ is sometimes called the {\sl variable speed random walk (VSRW)}. However, for the purpose of this paper it would also be possible to consider the {\sl constant speed random walk (CSRW)} with total jump rates normalized to one (cf.\ Remark~\ref{rem:csrw} below).

On $(\Omega, \mathcal{F},\mathbb{P})$ we define a $d+1$ parameter group of transformations $(\tau_{t,x})_{(t,x)\in \mathbb{R} \times \bZ^d}$ by
\[
\tau_{t,x}: \, \Omega\rightarrow\Omega \quad (\mu_e(s))_{s\in \mathbb{R}, e\in E_d}\mapsto (\mu_{x+e}(t+s))_{s\in \mathbb{R}, e\in E_d}
\]
so that obviously  $\tau_{s+t, x+y}=\tau_{s,x}\circ \tau_{t,y}$.  Notice that 
\begin{align}\label{eq:shiftedmu}
p^{\tau_{h,z}\om}(s,x;t,y)=p^\om(s+h,x+z;t+h, y+z), \qquad \mu_{xy}^{\tau_{h,z}\om}(t)=\mu^\om_{x+z,y+z}(t+h).
\end{align}
 
We are interested in the $\bP$ almost sure or quenched  long range behavior,  in particular in obtaining a
quenched functional limit theorem (QFCLT) or invariance principle
for the process $X$ starting in $0$ at time $0$. To that aim we need to state some assumptions on the environment measure $\bP$.

\begin{assumption}[Ergodicity] \label{ass:ergodic}
$\tau_{t,x}(A)\in \mathcal{F}$ for all $A\in \mathcal{F}$,
and the measure $\mathbb{P}$ is invariant and ergodic w.r.t.\ $(\tau_{t,x})$, i.e.\ $\mathbb{P}[A]\in\{0,1\}$ for any event $A$ such that $\tau_{t,x}(A)=A$ for all $t\in \bR$ and $x\in \bZ^d$.
\end{assumption}

\begin{assumption} [Stochastic Continuity]\label{ass:stoch_cont}
For any $\delta>0$ and $f\in L^2(\bP)$ we have 
\[
\lim_{h\to 0} \bP[ |f(\tau_{h,0}\om)-f(\om)|\geq \delta ]=0. 
\]
\end{assumption}

Thanks to Assumption \ref{ass:ergodic} and \ref{ass:stoch_cont}  the family of operators $(T_{t})_{t\in\bR}$ acting on $L^2(\bP)$, defined by $T_{t}f=f\circ \tau_{t,0}$, forms a strongly continuous group of unitary operators.  Its $L^2(\bP)$-generator will be denoted by $D_t: \, \mathcal{D}(D_t)\rightarrow L^2(\bP)$, defined by
\[
D_tf(\om)=\frac{\partial}{\partial t} T_tf_{|t=0}(\om)=\frac{\partial}{\partial t}_{|t=0} f(\tau_{t,0}\om).
\]
By Corollary 1.1.6 in \cite{EK} the generator is closed and densely defined. Note that $D_t$ is an anti-selfadjoint operator in $L^2(\bP)$, i.e.\
\[
\langle D_tf,g\rangle_{\bP}=-\langle f, D_t g\rangle_{\bP}, \qquad f,g\in \mathcal{D}(D_t),
\]
in particular
\begin{align} \label{eq:antisymD}
\langle D_tf,f\rangle_{\bP}=0,  \qquad f\in \mathcal{D}(D_t).
\end{align}

\begin{assumption} [Ellipticity] \label{ass:ellipticity}
There exist positive constants $C_l$ and $C_u$ such that
\begin{align} \label{ass_elliptic}
\mathbb{P}\big[  C_l \leq \mu_e(t) \leq C_u, \,\forall e\in E_d, t\in\mathbb{R} \big]=1.
\end{align}
\end{assumption}

We recall that under Assumption \ref{ass:ellipticity} the following heat kernel estimates have been proven in \cite{DD} (see also \cite[Appendix B]{GOS} for similar bounds).

\begin{proposition} \label{prop:hke}
There exist constants $c_1,\ldots,  c_5$ such that for $\bP$-a.e.\ $\om$ and for every $t\geq s\geq 0$ the following  holds:
\begin{enumerate}
\item[ i)] If $x,y\in \bZ^d$ and $D=|x-y|\leq c_1 (t-s)$, then
\begin{align*}
p^\om(s,x;t,y)\leq \frac{c_2} {(t-s)^{d/2}} \exp(-c_3 D^2/(t-s))\quad  \text{(Gaussian regime)}.
\end{align*}
\item [ii)] If $x,y\in \bZ^d$ and $D=|x-y|\geq c_1 (t-s)$, then
\begin{align*}
p^\om(s,x;t,y)\leq \frac{c_4}{1\vee (t-s)^{d/2}}  \exp(-c_5 D(1+ \log(D /(t-s))) \quad  \text{(Poisson regime)}.
\end{align*}
\end{enumerate}
\end{proposition}

Our first result is the following averaged or annealed FCLT. Let $\bP\otimes P_{s,x}^\om$ be the joint law of the environment and the random walk, and the annealed law is defined to be the marginal $\bP^*_{s,x}=\int_\Omega P_{s,x}^\om \, d\bP(\om)$. Further, let
\begin{equation*}
  X^{(\eps)}_t = \eps X_{t/\eps^2},  \q t \ge 0.
\end{equation*}

\begin{theorem}\label{annealed-ip} Let $d\geq 1$ and suppose that Assumptions \ref{ass:ergodic}-\ref{ass:ellipticity} hold.  Then, 
the law  of $X^{(\eps)}$ converges under $\bP^*_{0,0}$
to the law of a Brownian motion on $\bR^d$ with a deterministic non-degenerate covariance matrix
$\Sigma$.
\end{theorem}

To prove a QFCLT we will need some mixing assumptions on the environment. We denote by $B(\Omega)$ the set of bounded and measurable functions on $\Omega$ and 
$C^1_{b,\mathrm{loc}}(\hat \Omega)$ the set of differentiable functions on $\hat \Omega=[0,\infty)^{E_d}$ with bounded derivatives depending only on a finite number of variables.

\begin{assumption}[Time-mixing of the environment]\label{ass:time_mixing}
There exists $p_1>1$ such that for every $m\in\bN$   the following holds: For each $\varphi,\psi\in B(\Omega)$ of the form $\varphi(\om)=\tilde\varphi(\om(t_1))$ and $\psi(\om)=\tilde\psi(\om(t_2))$ with $|t_1-t_2|\geq 1$ for some $\tilde\varphi, \tilde\psi\in C^1_{b,\mathrm{loc}}(\hat\Omega)$ depending on $m$ variables we have
$$
\left| \bE[\varphi \psi ] - \bE[\varphi] \bE[\psi]  \right| \leq c_m |t_1-t_2|^{-p_1} \|\varphi\|_{L^\infty(\bP)} \|\psi\|_{L^\infty(\bP)}.
$$
\end{assumption}

\begin{assumption}[Space-mixing of the environment]\label{ass:space_mixing}
Let $d\geq 3$. There exists  $p_2>2d/(d-2)$ such that for every $m\in\bN$ and for every $x\in \bZ^d$  the following holds: For each $\varphi,\psi\in B(\Omega)$ of the form $\varphi(\om)=\tilde\varphi(\om(t_0))$ and $\psi(\om)=\tilde\psi(\om(t_0))$  for some $\tilde\varphi, \tilde\psi\in C^1_{b,\mathrm{loc}}(\hat\Omega)$ depending on $m$ variables we have
$$
\left| \bE[\varphi(\om) \psi(\tau_{0,x}\om) ] - \bE[\varphi] \bE[\psi]  \right| \leq c_m |x|^{-p_2} \|\varphi\|_{L^\infty(\bP)} \|\psi\|_{L^\infty(\bP)}.
$$
\end{assumption}

We are now ready to state the following QFCLT as our main result. 
\begin{theorem}\label{main-ip} Let $d\geq 3$ and suppose that Assumptions \ref{ass:ergodic}-\ref{ass:space_mixing} hold.  Then, 
$\bP$-a.s.\ 
$X^{(\eps)}$ converges (under $P_{0,0}^\om$)
in law to a Brownian motion on $\bR^d$ with a deterministic non-degenerate covariance matrix
$\Sigma$.
\end{theorem}

Notice that Theorem \ref{main-ip} only covers the transient lattice dimensions $d\geq 3$. In order to get an invariance principle for $X$ also in dimensions $d\leq 2$, we need to modify the mixing assumptions as follows.

\setcounter{assumptionpr}{3}
\begin{assumptionpr} \label{ass:time_mixing2}
Assumption \ref{ass:time_mixing} holds with $p_1>d+1$ if $d\geq 2$ and $p_1>4$ if $d=1$.
\end{assumptionpr}

\begin{assumptionpr} \label{ass:space_mixing2}
There exists $p_2>1$ such that for every $m\in\bN$ and for every $L>0$  the following holds: For each  $\varphi,\psi\in B(\Omega)$ of the form $\varphi(\om)=\tilde\varphi(\om(t_1))$ and $\psi(\om)=\tilde\psi(\om(t_2))$, where $|t_1-t_2|\leq L$ and $\tilde\varphi, \tilde\psi\in C^1_{b,\mathrm{loc}}(\hat\Omega)$ depend on variables contained in two subsets $A_\varphi$ and $A_\psi$ of $\bZ^d$ with diameter at most $m$ and $\dist (A_\varphi, A_\psi)\geq L$, 
$$
\left| \bE[\varphi \psi ] - \bE[\varphi] \bE[\psi]  \right| \leq c_m L^{-p_2} \|\varphi\|_{L^\infty(\bP)} \|\psi\|_{L^\infty(\bP)}.
$$
\end{assumptionpr}

\begin{theorem}\label{main-ip2} Let $d\geq 1$ and suppose that Assumptions \ref{ass:ergodic}-\ref{ass:ellipticity}, \ref{ass:time_mixing2} and \ref{ass:space_mixing2} hold. 
 Then, 
$\bP$-a.s. 
$X^{(\eps)}$ converges (under $P_{0,0}^\om$)
in law to a Brownian motion on $\bR^d$ with a deterministic non-degenerate covariance matrix
$\Sigma$.
\end{theorem}

\begin{remark} \label{rem:csrw}
One can also consider the time-inhomogeneous {\sl constant speed random walk or CSRW} $Y=(Y_t, t\in \bR, P_{s,x}^\omega, (s,x)\in \bR \times\bZ^d)$ with generator given by:
\begin{align*}
 \sL^{Y}_t f(x) = \sum_{y\sim x} \frac{\mu_{xy}(t)}{\mu_x(t)} (f(y)-f(x)).
\end{align*}
In contrast to the VSRW $X$, whose waiting time at any site $x\in \bZ^d$ depends on $x$, the CSRW waits at each site an exponential time with mean one. Since the CSRW is a time change of the VSRW, an invariance principle for $Y$ follows from an  invariance principle for $X$  by the same arguments as in \cite[Section 6.2]{ABDH}. In this case the limiting object is a Brownian motion in $\mathbb{R}^d$ with covariance matrix $\Sigma_C=(1/\bE \mu_0(0)) \Sigma_V$, where $\Sigma_V$ denotes the covariance matrix of the limiting Brownian motion in the invariance principle for $X$.
\end{remark}

Next we state some consequences of our results, which follow from arguments in \cite{BH} by combining the invariance principle for $X$ and the Gaussian bound for the heat kernel. First, we have a local  limit theorem for the heat kernel. 
 Write  
\begin{align} \label{def_k}
 k_t(x) =k_t^{(\Sigma)}(x)= \frac 1 {\sqrt{(2\pi t)^d \det \Sigma } } \exp(- x\cdot \Sigma^{-1}x /2t)
 \end{align}
for the Gaussian heat kernel with diffusion matrix $\Sigma$.

\begin{theorem}\label{thm:llt}  Let $T>0$. 
For $x \in \bR^d$ write 
$\lfloor x \rfloor =(  \lfloor x_1 \rfloor, \dots \lfloor x_d \rfloor)$.
\begin{enumerate}
\item[i)] Suppose that Assumptions \ref{ass:ergodic}-\ref{ass:ellipticity} hold. Then,
\begin{align*}
 \lim_{n\to\infty}  
\sup_{x\in \bR^d} \sup_{t \ge T} 
 \Big|n^{d/2}\bE[ p^\om(0,0; nt, \lfloor n^{1/2}x \rfloor)] -  k_t(x)\Big| =0.
\end{align*}
\item[ii)] Under the assumptions of Theorem  \ref{main-ip} or Theorem \ref{main-ip2} we have
\begin{align*} 
  \lim_{n\to\infty}  
\sup_{x\in \bR^d} \sup_{t \ge T} 
 \Big|n^{d/2}p^\om(0,0; nt, \lfloor n^{1/2}x \rfloor)-  k_t(x)\Big| =0, \q 
\hbox{$\bP$-a.s.} 
\end{align*}
\end{enumerate}
\end{theorem}

\proof Given the annealed or quenched invariance principle and the
heat kernel bounds in Proposition~\ref{prop:hke} this can be proven as in Section 4 of \cite{BH}.  
\qed

When $d\geq 3$ the calculations in Section~6 of \cite{BH} then give the following bound on the Green kernel $g^\om(x,y)$ defined by
\[
 g^\om (x,y)=\int_0^\infty p^\om (0,x;t,y) \, dt.
\]

\begin{theorem}
 Let $d\geq 3$ and suppose that the assumptions of Theorem~\ref{main-ip} or Theorem~\ref{main-ip2} hold.
\begin{enumerate}
 \item[i)] There exist constants $c_1$ and $c_2$ such that for $x\not= y$
\[
 \frac{c_1}{|x-y|^{d-2}} \leq g^\om (x,y) \leq \frac{c_2}{|x-y|^{d-2}}.
\]
\item[ii)] Let $C=\Gamma(\frac d 2 -1)/2 \pi^{d/2} \det \Sigma$. For any $\varepsilon >0$ there exists $M=M(\varepsilon, \om)$ with $\bP[M<\infty]=1$ such that
\begin{align*}
\frac{(1-\varepsilon) C}{|x|^{d-2}} \leq g^\om (0,x) \leq \frac{(1+\varepsilon)C}{|x|^{d-2}} \qquad \text{for $|x|>M(\om)$.}
\end{align*}

\item[iii)] We have, $\bP$-a.s.,
\[
 \lim_{|x|\to\infty} |x|^{2-d} g^\om(0,x)= \lim_{|x|\to\infty} |x|^{2-d} \bE [g^\om(0,x)]=C.
\]
\end{enumerate}
\end{theorem}

In the case of static conductances, quenched invariance principles for the random conductance model have been proven by a number of different authors under various restrictions on the law of the conductances, see \cite{SS,BP,Ma,BD}. Recently, these results have been unified in \cite{ABDH}, where a QFCLT has been obtained for the RCM with general nonnegative i.i.d.\ conductances.  We also refer the reader to \cite{Bi} for a recent survey on this topic.

On the other hand, to our knowledge the present paper is the first one proving an invariance principle for the RCM with a time-dynamic environment.  
However, quenched invariance principles have been proven for several other discrete-time random walks in a dynamic random environment. In \cite{BMP1} a QFCLT is obtained for random walks in space-time product environments  by using Fourier-analytic methods. This result has been improved in \cite{BMP2} to environments satisfying an exponential spatial mixing assumption and in \cite{BZ} to Markovian environments by using more probabilistic techniques.
Another very successful approach is the well-established Kipnis-Varadhan technique based on the process of the environment as seen from the particle. 
In \cite{RS} this approach has been used to get a QFCLT for the random walk in space-time product environments. Moreover, it has been applied in \cite{DL2} to random walks in a dynamic enviroment, which forms a Gibbsian Markov chain in time with spatial mixing, and in \cite{JR} to random walks on $\mathbb{R}^d$, where the environment is i.i.d.\ in time and polynomially mixing in space. 
Recently, a general class of random walks in an ergodic Markovian environment satisfying some coupling conditions has been studied in \cite{RV}.

Also in this paper we will follow the approach in \cite{RS}, so we use the process of the environment as seen from the particle and the method of the 'corrector', that is we decompose the random walk $X$ into a martingale and a time-dependent corrector function.   Due to the time-inhomogeneity  and the resulting lack of reversibility we need to apply the adaptions of the Kipnis-Varadhan method to non-reversible situations in \cite{MW} and \cite{KLO}. In particular, in order to construct the corrector we show that the generator of the environment seen from the particle is a perturbation of a normal operator in the sense of \cite[Section 2.7.5]{KLO}. This is done in Section~\ref{sec:corr}.  As a byproduct this will already imply the annealed FCLT in Theorem~\ref{annealed-ip}.

Once the corrector is constructed, the QFCLT for the martingale part is standard, so it remains to control the corrector. To that aim we still follow \cite{RS} and apply the theory of 'fractional coboundaries' of Derriennic and Lin in \cite{DL}. The main step in this approach is to establish a subdiffusive bound on the corrector (see Proposition~\ref{prop:asymp_corr} below), which is done in Section~\ref{sec:asymp_corr}. To obtain this bound we establish so-called two-walk estimates, i.e.\ we consider the difference of two independent copies of $X$ evolving in the same fixed environment $\om$ (cf.\ e.g.\ \cite{JR} or Appendix A in \cite{RS2}). In $d\geq 3$, following \cite{Mou}  we show that the variance decay of the environment viewed from the particle is strong enough for our purposes by using the mixing assumption \ref{ass:time_mixing} and \ref{ass:space_mixing}  (see Lemma \ref{lem:asymp_Pk}). In the recurrent lattice dimensions $d\leq 2$ the estimate for the variance decay is not good enough, so we give a different argument here involving the modified mixing asumptions  \ref{ass:time_mixing2} and \ref{ass:space_mixing2}.

In Section \ref{sec:pf_main} we prove the main result, i.e.\ we state a tightness result, which is a direct consequence from the heat kernel bounds in Proposition~\ref{prop:hke}, and show the QFCLT for the martingale part. To control the corrector we apply the results in \cite{DL}, which are stated in the discrete-time setting. Since it is not clear to us, how to apply them directly in the continuous-time setting, we first prove the QFCLT for the discretized process as in \cite{BD}. More precisely, we 
define $\wh X_n = X_n$, $n \in \bN$, and consider the process
\begin{equation*}
 \wh X^{(\eps)}_t =  \eps \wh X_{ \lfloor t/\eps^2 \rfloor}.
\end{equation*}
We can control $\sup_{t \le T} | X^{(\eps)}_t - \wh X^{(\eps)}_t|$
-- see Lemma \ref{disc-approx} --
so an invariance principle for $X^{(\eps)}$ will follow 
from one  for $\wh X^{(\eps)}$. 

Finally, in Section \ref{sec:interface}  we  point out a link to stochastic interface models (see \cite{F}). Namely, a local limit theorem for the RCM with dynamic conductances can be used to obtain scaling limits for the space-time covariation of the Ginzburg-Landau interface model via Helffer-Sj\"ostrand representation. 

Throughout the paper we write $c$ to denote a positive constant which may change
on each appearance. Constants denoted $c_i$ will be the same through
each argument.

\bigskip
\textbf{Acknowledgement.} I thank Martin Barlow, Jean-Dominique Deuschel and Martin Slowik for helpful discussions and useful comments.

\section{Construction of the Corrector} \label{sec:corr}
Throughout this section we suppose that Assumptions \ref{ass:ergodic}-\ref{ass:ellipticity} hold. We define the process of the environment seen from the particle by
\[
\eta_t(\om)=\tau_{t,X_t}\om, \qquad \om \in \Omega, \, t\geq 0.
\]

\begin{proposition} \label{prop:semigroup_eta}
\begin{enumerate}
\item[i)] The process $(\eta_t)_{t\geq 0}$ is Markovian with transition semigroup
\[
P_tf(\om)=\sum_{y \in\bZ^d}  p^\om(0,0;t,y) f(\tau_{t,y}\om) \qquad \mbox{for all $f\in B(\Omega)$}.
\]
The semigroup $(P_t)$ extends uniquely to a strongly continuous semigroup of contractions on $L^2(\bP)$, whose generator $L: \mathcal{D}(L)\rightarrow L^2(\bP)$ is given by
\begin{align*} 
Lf(\om)=D_tf(\om)+\sum_{y\sim 0} \mu^\om_{0y}(0) (f(\tau_{0,y}\om)-f(\om))
\end{align*}
with domain $\mathcal{D}(L)=\mathcal{D}(D_t)$.
\item [ii)] The measure $\bP$ is invariant and ergodic for $\eta$.
\end{enumerate}
\end{proposition}
\proof
i) The Markov property as well as the representation of the semigroup follow from \eqref{eq:shiftedmu} by similar arguments as in Lemma 3.1 in \cite{KLO}. For every bounded $f\in\mathcal{D}(D_t)$ we have
\begin{align*}
\frac  {P_tf(\om)-f(\om)} t=\sum_{y\in\bZ^d} \frac{p^\om(0,0;t,y)}{t} (f(\tau_{0,y}\om)-f(\om))+\sum_{y\in\bZ^d} p^\om(0,0;t,y) \frac{f(\tau_{t,y}\om)-f(\tau_{0,y}\om)}{t}.
\end{align*}
Taking limits for $t\downarrow 0$, using the fact that $ p^\om(0,0;t,y)\to \delta_{0y}$, we obtain the formula for $Lf$. Obviously, the operators $L$ and $D_t$ have the same domain.

ii) Let $f\in\mathcal{D}(L)$. Since the operator $D_t$ is anti-selfadjoint we have $\langle D_t f\rangle_{\bP}=0$. Hence,
\begin{align*}
\langle L f\rangle_{\bP}&=\sum_{y\in\bZ^d} \langle \mu_{0y}^\om(0) f(\tau_{0,y}\om) \rangle_{\bP}-\langle \mu_{0y}^\om(0) f(\om) \rangle_{\bP}
=\sum_{y\in\bZ^d} \langle \mu_{0y}^{\tau_{0,-y}\om}(0) f(\om) \rangle_{\bP}-\langle \mu_{0y}^\om(0) f(\om) \rangle_{\bP} \\
&=\sum_{y\in\bZ^d} \langle \mu_{0,-y}^\om(0) f(\om) \rangle_{\bP}-\langle \mu_{0y}^\om(0) f(\om) \rangle_{\bP}=0,
\end{align*}
where we have used the invariance of $\bP$ w.r.t.\ $\tau_{t,x}$ and \eqref{eq:shiftedmu}. Thus, $\bP$ is an invariant measure for $\eta$. To prove that $\bP$ is also ergodic, let now $A\in \mathcal{F}$ with $P_t \indicator_A=\indicator_A$. Then,
\begin{align*}
0=\indicator_{A^c}(\om) \cdot P_t \indicator_A(\om)=\sum_{y\in\bZ^d}  \indicator_{A^c}(\om) p^\om(0,0;t,y) \indicator_A(\tau_{t,y}\om).
\end{align*}
Since for all $t>0$ and $y\in \bZ^d$ there is a stricly positive lower bound for $p^\om(0,0;t,y)$ independent of $\om$ (see Proposition~4.3 in \cite{DD}) we get
\[
\indicator_{A^c}(\om) \cdot \indicator_A(\tau_{t,y}\om)=0.
\]
Thus, the set $A$ is invariant under $\tau_{t,x}$. Since $\bP$ is ergodic w.r.t.\ $\tau_{t,x}$ we conclude that $A$ is $\bP$-trivial and the claim follows.
\qed

\begin{lemma} \label{lem:DF_eta}
For $f\in \mathcal{D}(L)$,
\[
\langle f,(-L)f\rangle_{\bP}=\tfrac 1 2 \sum_{y\in\bZ^d} \bE \left[ \mu_{0y}^\om(0) (f(\tau_{0,y}\om) -f(\om))^2 \right].
\]
\end{lemma}
\proof
Recall that $\langle f, D_tf \rangle_{\bP}=0$. Therefore,
\begin{align*}
&\langle f,(-L)f \rangle_{\bP} 
=-\sum_{y\in\bZ^d} \bE \left[ f(\om)  \mu_{0y}^\om (0) (f(\tau_{0,y}\om)-f(\om)) \right] \\
=& -\tfrac 1 2 \sum_{y\in\bZ^d} \bE \left[ f(\om)  \mu_{0y}^\om (0) (f(\tau_{0,y}\om)-f(\om)) \right]-\tfrac 1 2 \sum_{y\in\bZ^d} \bE \left[ f(\om)  \mu_{0,-y}^\om (0) (f(\tau_{0,-y}\om)-f(\om)) \right] \\
=& -\tfrac 1 2 \sum_{y\in\bZ^d} \bE \left[ f(\om)  \mu_{0y}^\om (0) (f(\tau_{0,y}\om)-f(\om))\right]-\tfrac 1 2 \sum_{y\in\bZ^d} \bE \left[ f(\tau_{0,y}\om)  \mu_{0,-y}^{\tau_{0,y}\om} (0) (f(\om)-f(\tau_{0,y}\om)) \right] \\
=&\tfrac 1 2 \sum_{y\in\bZ^d} \bE \left[ \mu_{0y}^\om(0) (f(\tau_{0,y}\om) -f(\om))^2 \right],
\end{align*}
where we have used again the invariance of $\bP$ w.r.t.\ $\tau_{t,x}$ and \eqref{eq:shiftedmu}. 
\qed

Let $P_t^*$ and $L^*$ denote the $L^2(\bP)$-adjoint operators of $P_t$ and $L$, respectively.
\begin{proposition} \label{prop:adjoint}
We have
\[
P_t^*f(\om)=\sum_{y\in\bZ^d} \hat p^\om(0,0;t,y) f(\tau_{-t,y}\om), \qquad \mbox{$f\in L^2(\bP)$},
\]
with $\hat p^\om (s,x;t,y):=p^\om(-t,y;-s,x)$ and for $f\in \mathcal{D}(L)$
\[
L^*f(\om)=-D_tf(\om)+\sum_{y\sim 0} \mu^\om_{0y}(0) (f(\tau_{0,y}\om)-f(\om)).
\]
\end{proposition}
\proof 
Using \eqref{eq:shiftedmu} we compute the adjoint of $P_t$ as
\begin{align*}
\langle P_t f, g\rangle_{\bP}&= \sum_{y\in\bZ^d} \bE \left[ p^\om(0,0;t,y) f(\tau_{t,y}\om) g(\om)\right]= \sum_{y\in\bZ^d} \bE \left[ p^{\tau_{-t,-y}\om}(0,0;t,y) f(\om) g(\tau_{-t,-y}\om)\right] \\
&=\sum_{y\in\bZ^d} \bE \left[ p^{\om}(-t,-y;0,0) f(\om) g(\tau_{-t,-y}\om)\right]=\sum_{y\in\bZ^d} \bE \left[ \hat p^{\om}(0,0;t,y) g(\tau_{-t,y}\om)  f(\om) \right],
\end{align*}
and the representation for $P_t^*$ follows. To compute $L^*$ we use a similar procedure as in Lemma~\ref{lem:DF_eta} and get
\begin{align*}
\langle Lf,g \rangle_{\bP}&=\langle D_t f,g \rangle_{\bP}+\sum_{y\in\bZ^d} \bE\left[ \mu_{0y}^\om(0) (f(\tau_{0,y}\om)-f(\om)) g(\om) \right] \\
&=-\langle f,D_tg \rangle_{\bP}-\tfrac 1 2 \sum_{y\in\bZ^d} \bE\left[ \mu_{0y}^\om(0) (f(\tau_{0,y}\om)-f(\om)) (g(\tau_{0,y}\om)-g(\om))  \right] \\
&=-\langle f,D_tg \rangle_{\bP}+\sum_{y\in\bZ^d} \bE\left[ \mu_{0y}^\om(0) (g(\tau_{0,y}\om)-g(\om)) f(\om) \right],
\end{align*}
which gives the claim.
\qed

Next we introduce the Hilbert spaces $\mathcal{H}_1$ and $\mathcal{H}_{-1}$. Let $\mathcal{C}$ be a common core of the operators $L$ and $L^*$. On $\mathcal{C}$ we define the seminorm 
\[
\| f\|^2_{\mathcal{H}_1}=\langle f,(-L)f\rangle_{\bP}, \qquad f\in \mathcal{C}.
\]
Let $\mathcal{H}_1$ be the completion of $\mathcal{C}$ (or more precisely the completion of equivalence classes of elements in $\mathcal{C}$ w.r.t.\ the equivalence relation $f\sim g$ if $\|f-g\|_{\mathcal{H}_1}=0$) w.r.t.\ $\|.\|_{\mathcal{H}_1}$. Then, $\mathcal{H}_1$ is a Hilbert space  with inner product $\langle .,. \rangle_{\mathcal{H}_1}$ given by polarization:
\[
 \langle f,g \rangle_{\mathcal{H}_1}=\frac 1 4 \left( \| f+g \|^2_{\mathcal{H}_1}- \| f-g\|^2_{\mathcal{H}_1} \right).
\]
Associated with $\mathcal{H}_1$ we define the dual space $\mathcal{H}_{-1}$ as follows. For $f\in L^2(\bP)$ let
\[
 \| f\|^2_{\mathcal{H}_{-1}}=\sup_{g\in \mathcal{C}}\left( 2 \langle f ,g \rangle_{\mathbb{P}} - \| g\|^2_{\mathcal{H}_1} \right).
\]
The Hilbert space $\mathcal{H}_{-1}$ is then defined as the $\| .\|_{\mathcal{H}_{-1}}$-completion of (equivalence classes of) elements in $\mathcal{C}$ with finite $\| .\|_{\mathcal{H}_{-1}}$-norm. As before the inner product  $\langle .,. \rangle_{\mathcal{H}_{-1}}$ is defined through polarization.
We refer to Section 2.2 in \cite{KLO} for more details.

Next we define the local drift
\[
V_j(\om)=\sum_{y\sim 0} \mu_{0y}^\om (0) y^j =\sL^\om_0 f_j(0), \qquad j=1,\ldots,d,
\]
where $f_j(x)=x^j$, $x^j$ and $y^j$ denoting the $j$-th component of $x$ and $y$. Since $\mu^\om_{0y}(0)=0$ unless $y\sim 0$, we have $V_j(\om)=\mu^\om_{0,e_j}(0)-\mu^\om_{0,-e_j}(0)$.

\begin{lemma} \label{lem:V_H-1}
For every $j=1,\ldots,d$,  $V_j\in L^2(\bP)\cap \mathcal{H}_{-1}$.
\end{lemma}
\proof
It suffices to show that 
\begin{align} \label{eq:estH-1}
\big| \langle V_j,f\rangle_{\bP}\big|^2 \leq c \langle f, (-L)f \rangle_{\bP} \qquad \text{for all $f\in\mathcal{H}_1$}
\end{align}
 (cf.\ equation (2.12) in \cite{KLO})). By definition of $V_j$ we have
\begin{align*}
\langle V_j,f\rangle_{\bP}&=\bE \left[ \mu^\om_{0,e_j}(0) f(\om) \right]- \bE \left[\mu^\om_{0,-e_j}(0) f(\om)\right]=\bE \left[ \mu^\om_{0,e_j}(0) f(\om) \right]- \bE \left[ \mu^\om_{0,e_j}(0) f(\tau_{0,e_j}\om) \right] \\
&=-\bE \left[ \mu^\om_{0,e_j}(0) (f(\tau_{0,e_j}\om)-f(\om)) \right].
\end{align*}
Hence, using Cauchy Schwarz and Lemma \ref{lem:DF_eta}
\begin{align*}
\big| \langle V_j,f\rangle_{\bP}\big|^2 &\leq  \bE[ \mu^\om_{0,e_j}(0)] \, \bE[ \mu^\om_{0,e_j}(0) (f(\tau_{0,e_j}\om)-f(\om))^2]
 \leq C_u \sum_{y\in\bZ^d} \bE \left[ \mu^\om_{0y}(0) (f(\tau_{0,y}\om)-f(\om))^2 \right] \\
 &=2 C_u  \langle f, (-L)f \rangle_{\bP},
\end{align*}
and we obtain \eqref{eq:estH-1}.
\qed

For $\lambda>0$, we  consider for each $j$ the solution $u_\lambda^{j}$ of the resolvent equation 
\begin{equation} \label{eq:res}
(\lambda-L) u_\lambda^{j}=V_j.
\end{equation}

\begin{proposition} \label{prop:conv_u}
For every $j=1,\ldots,d$, there exists $u^{j} \in \mathcal{H}_1$ such that
\begin{align*}
\lim_{\lambda\to 0} \lambda \| u_\lambda^{j} \|^2_{L^2(\bP)}=0 \quad \text{and} \quad 
\lim_{\lambda\to 0}u_\lambda^{j}=u^{j} \, \text{strongly in $\mathcal{H}_1$}.
\end{align*}
\end{proposition}

The proof of Proposition~\ref{prop:conv_u} will be based on the following statement proven in \cite{KLO}.
\begin{proposition} \label{prop:pert_op}
Suppose we have the decomposition  $L=L^0+B$ of the operator $L$ such that 
\begin{enumerate}
\item[i)] The operator $L^0$ is normal, i.e.\ $L^0 (L^0)^*=(L^0)^*L^0$.
\item[ii)] The Dirichlet forms of $L$ and $L^0$ are equivalent, i.e.\ there exist positive constants $c_1$ and $c_2$ such that
\[
c_1 \langle f, (-L)f \rangle_{\bP}\leq \langle f, (-L^0)f \rangle_{\bP}\leq c_2 \langle f, (-L)f \rangle_{\bP}, \qquad \text{for all $f\in \mathcal{D}(L)$}.
\]
\item[iii)] $B$ satisfies a sector condition w.r.t.\ $L^0$, i.e.\ there exists a positive constant $c$ such that
\begin{align*}
\langle f,Bg\rangle^2_{\bP} \leq c  \langle f, (-L^0)f\rangle_{\bP} \langle g, (-L^0)g\rangle_{\bP}, \qquad f,g\in \mathcal{D}(L).
\end{align*}
\end{enumerate}
Then, for any fixed $V\in L^2(\bP)\cap \mathcal{H}_{-1}$ the solution $f_\lambda$ of the resolvent equation
$(\lambda-L)f_\lambda=V$ satisfies
\begin{align*}
\lim_{\lambda\to 0} \lambda \| f_\lambda\|^2_{L^2(\bP)}=0 \quad \text{and} \quad 
\lim_{\lambda\to 0}f_\lambda=f \, \text{strongly in $\mathcal{H}_1$},
\end{align*}
for some  $f \in \mathcal{H}_1$.
\end{proposition}
\proof By Proposition 2.25 in \cite{KLO} the assumptions imply that
\[
\sup_{0<\lambda \leq 1} \|Lf_\lambda \|_{\mathcal{H}_{-1}}<\infty.
\]
The claim follows then from Lemma 2.16 in \cite{KLO}.
\qed

\textbf{Proof of Proposition~\ref{prop:conv_u}.}
We decompose the operator $L=L^0+B$ with 
\begin{align*}
L^0 f&:=D_t f+\sum_{y\sim 0} C_l (f(\tau_{0,y}\om)-f(\om)), \qquad f\in \mathcal{D}(L), \intertext{and} 
Bf&:= \sum_{y\sim 0} (\mu^\om_{0y}(0)-C_l) (f(\tau_{0,y}\om)-f(\om)), \qquad f\in \mathcal{D}(L).
\end{align*}
A similar calculation as in the proof of Lemma \ref{lem:DF_eta} shows that
\begin{align} \label{eq:DFL0}
\langle f, (-L^0)f \rangle_{\bP}&=\tfrac 1 2 \sum_{y\sim 0 } \bE \left[ C_l (f(\tau_{0,y}\om)-f(\om))^2 \right],\\ \label{eq:DFB}
\langle f, (-B)g \rangle_{\bP}&=\tfrac 1 2 \sum_{y\sim 0 }  \bE \left[ (\mu^\om_{0y}(0)-C_l) (f(\tau_{0,y}\om)-f(\om)) (g(\tau_{0,y}\om)-g(\om))\right].
\end{align}
The claim will follow from Proposition \ref{prop:pert_op} and Lemma \ref{lem:V_H-1} once we have verified conditions i)-iii) in Proposition~\ref{prop:pert_op}.
 To show i), note that the closure of $L^0$ is the generator of a semigroup $(P^0_t)$ that corresponds to a process seen from the particle associated with a simple random walk  on $\bZ^d$ with constant jump rates $C_l$. In particular, the associated process is time-homogeneous, i.e.\ the corresponding transition probabilities satisfy $p_0^\om(s,x;t,y)=p_0^\om(t-s,x,y)$ and $\hat p_0^\om(s,x;t,y)=\hat p_0^\om(t-s,x,y)$, where $p_0^\om(t,x,y)=p_0^\om(0,x;t,y)$ and $\hat p_0^\om(t,x,y)=\hat p_0^\om(0,x;t,y)$. Since this random walk is obviously reversible w.r.t.\ the counting measure, we have $p_0^\om(t,x,y)=\hat p_0^\om(t,x,y)$. Then, since we have similar representations for $P^0_t$ and $(P^0_t)^*$ as for the semigroups in Proposition \ref{prop:semigroup_eta} and Proposition \ref{prop:adjoint}, we get
\[
(P^0_t)^*P^0_t=P^0_t (P^0_t)^*, \qquad t\geq 0,
\]
which implies that the closure of $L^0$ is normal (see Theorem 13.37 in \cite{Ru}).

Condition ii) is immediate from Lemma \ref{lem:DF_eta},  \eqref{eq:DFL0} and the ellipticity condition \eqref{ass_elliptic}. To prove iii) we use \eqref{eq:DFB}, Cauchy Schwarz and the ellipticity condition \eqref{ass_elliptic}, which gives
\begin{align*}
\langle f,Bg\rangle^2_{\bP} &\leq 
\tfrac 1 2 C_u^2 d \, \bE \bigl[ \sum_{y\sim 0 }  (f(\tau_{0,y}\om)-f(\om))^2 \bigr] \times  \bE \bigl[ \sum_{y\sim 0 }(g(\tau_{0,y}\om)-g(\om))^2 \bigr] \\
 &\leq  \frac{ C^2_u d}{2 C^2_l}  \langle f, (-L^0)f\rangle_{\bP} \langle g, (-L^0)g\rangle_{\bP},
\end{align*}
and the claim follows.
\qed

For abbreviation we write $u_\lambda=(u_\lambda^{1},\ldots,u_\lambda^{d})$ and  
\[
\chi_{\lambda}(t,x,\om):= u_\lambda \circ \tau_{t,x}-u_\lambda.
\]

\begin{proposition} \label{prop:constr_corr}
 For all non-negative $t\in \mathbb{Q}$ and $x\in \bZ^d$   the limit
\begin{align} \label{def_chi}
\lim_{\lambda'\to 0}\chi_{\lambda'}(t,x,\om) =:\chi(t,x,\om)
\end{align}
exists  along a subfamily $(\lambda')$ for $\bP$-a.e.\ $\om$. Moreover, the mapping $t\mapsto \chi(t,X_t,\om)$ can be extended to a right-continuous function on $[0,\infty)$ such that
\begin{align} \label{eq:decompX}
M_t=X_t+\chi(t,X_t,\om), \qquad t\geq 0,
\end{align}
is a $P_{0,0}^\om$-martingale.  
\end{proposition}

\proof
For every $j=1,\ldots,d$ and every $\lambda>0$ we have that for $\bP$-a.e.\ $\om$ the processes
\begin{align} \label{def_Nk}
N^{j,\lambda}_t=u_\lambda^{j}(\eta_t)-u_\lambda^{j}(\om)-\int_0^t Lu_\lambda^{j}(\eta_s) \, ds
\end{align}
and
\begin{align}\label{def_tildeM}
\tilde M^j_t=X^j_t-\int_0^t \mathcal{L}^\om_s f_j(X_s) \, ds
\end{align}
are both $P_{0,0}^\om$-martingales, where as before $f_j(x)=x^j$. Then, using the definition of $V_j$ and the fact that $u_\lambda^{j}$ solves the resolvent equation \eqref{eq:res} we get
\begin{align}\label{eq:mart_decomp}
X_t^j&=\tilde M^j_t+\int_0^t \mathcal{L}^\om_s f_j(X_s) \, ds=\tilde M^j_t+\int_0^t V_j(\eta_s) \, ds=\tilde M^j_t+\int_0^t (\lambda-L)u_\lambda^{j}(\eta_s) \, ds \nonumber \\
&=\tilde M^j_t+N^{j,\lambda}_t-\left(u_\lambda^{j}(\eta_t)-u_\lambda^{j}(\om)\right)+\lambda \int_0^t u_\lambda^{j}(\eta_s) \, ds.
\end{align}

In a first step we show that the martingale $N^{j,\lambda}_t$ converges in $L^2(\bP\otimes P_{0,0}^\om)$ as $\lambda \downarrow 0$ to a martingale $N^j_t$. To that aim it is enough to prove that $N^{j,\lambda}_t$ is a Cauchy sequence in $L^2(\bP\otimes P_{0,0}^\om)$. Since $\bP$ is an invariant measure for $\eta$ we use Lemma \ref{lem:DF_eta} to obtain
\begin{align*}
\bE E_{0,0}^\om \langle N^{j,\lambda}-N^{j,\lambda'}\rangle_t&= \int_0^t \bE E_{0,0}^\om \left[ L (u_\lambda^{j}-u_{\lambda'}^{j})^2-2 (u_\lambda^{j}-u_{\lambda'}^{j}) L(u_\lambda^{j}-u_{\lambda'}^{j}) \right](\eta_s) \, ds \\
&=2 t \langle (u_\lambda^{j}-u_{\lambda'}^{j}), (-L)(u_\lambda^{j}-u_{\lambda'}^{j})\rangle_{\bP}  \\
&=2t \| u_\lambda^{j}-u_{\lambda'}^{j} \|^2_{\mathcal{H}_1},
\end{align*}
which implies that $N^{j,\lambda}_t$ is a Cauchy sequence in $L^2(\bP\otimes P_{0,0}^\om)$ by Proposition \ref{prop:conv_u}. Thus, the martingale  $M_t^{j,\lambda}=\tilde M_t^j+N_t^{j,\lambda}$ converges to a martingale, whose right-continuous modification we denote by $M^j_t$. We define $M_t=(M_t^1,\ldots, M_t^d)$.

The next step is to show that the last term in \eqref{eq:mart_decomp} converges to zero in $L^2(\bP\otimes P_{0,0}^\om)$ as $\lambda \downarrow 0$. Since $V_j\in \mathcal{H}_{-1}$ we have that $\lim_\lambda \lambda u_\lambda^{j}=0$ in $L^2(\bP)$ (cf.\ equation (2.15) in \cite{KLO}). Thus, for every $j=1,\ldots,d$,
\begin{align*}
\big\| \lambda  \int_0^t u_\lambda^{j}(\eta_s) \, ds \big\|_{L^2(\bP\otimes P_{0,0}^\om)} \leq  \lambda \int_0^t \| u_\lambda^{j}(\eta_s)  \|_{L^2(\bP\otimes P_{0,0}^\om)}\, ds=t \, \lambda \| u_\lambda^{j}  \|_{L^2(\bP)} \rightarrow 0.
\end{align*}

Thus, by taking $L^2(\bP\otimes P_{0,0}^\om)$-limits in \eqref{eq:mart_decomp}  we get that $\chi_{\lambda}(t,X_t,\cdot)$ converges  in $L^2(\bP\otimes P_{0,0}^\om)$ as $\lambda\downarrow 0$ for every $t\geq 0$.  By a diagonal procedure we can extract a suitable subsequence $\lambda'$ such that for $\bP$-a.e.\ $\om$ we have that 
$\chi_{\lambda}(t,X_t,\om)$  has a limit in $L^2( P_{0,0}^\om)$ and $P_{0,0}^\om$-a.s.\ along  $\lambda'$ for all non-negative $t\in\mathbb{Q}$.  In particular, the limit is $\sigma(X_t)$-measurable and will therefore be denoted by $\chi(t,X_t,\om)$. Hence, 
\begin{align}\label{eq:Mpre}
X_t=M_t-\chi(t,X_t,\om).
\end{align}
Moreover, for $\bP$-a.e.\ $\om$,
\begin{align*}
\sum_{y\in\bZ^d} p^\om(0,0;t,y) \left| \chi_{\lambda}(t,y,\om) - \chi(t,y,\om) \right|^2=E^0_\om  \left| ( \chi_{\lambda}(t,X_t,\om) ) - \chi(t,X_t,\om)\right|^2 \rightarrow 0
\end{align*}
along $\lambda'$.
 Since $p^\om(0,0;t,y)>0$ for all $t>0$ and $x\in \bZ^d$, we conclude that for $\bP$-a.e.\ $\om$ the limit in \eqref{def_chi} exists for every non-negative $t\in \mathbb{Q}$ and every $y\in \mathbb{Z}^d$. Finally, using \eqref{eq:Mpre}  and the fact that $X_t$ and $M_t$ have right-continuous trajectories, we can extend $\chi(t,X_t,\om)$  to a right-continuous function on $[0,\infty)$  and \eqref{eq:decompX} follows. 
\qed

\begin{remark} \label{rem:defh}
Note that for all non-negative $s,t \in \mathbb{Q}$ and $x,y\in\bZ^d$,
\begin{align*}
u_\lambda(\tau_{t,y}\om)-u_\lambda(\tau_{s,x}\om)=&\left(u_\lambda(\tau_{t,y}\om)-u_\lambda(\om) \right)- \left(u_\lambda(\tau_{s,x}\om)-u_\lambda(\om)\right) \\
\rightarrow & \chi(t,y,\om)-\chi(s,x,\om)
\end{align*}
along the chosen subsequence for $\bP$-a.e.\ $\om$. The function $h_\lambda(\om_0,\om_1):=u_\lambda(\om_1)-u_\lambda(\om_0)$ on $\Omega \times \Omega$ converges in $L^2(\Omega\times\Omega, \bP\circ (\tau_{s,x},\tau_{t,y})^{-1})$ to a function $h$.  In particular, for $\bP$-a.e.\ $\om$,
\[
h(\tau_{s,x}\om,\tau_{t,y}\om)=\chi(t,y,\om)-\chi(s,x,\om).
\]
\end{remark}

\begin{corollary} \label{cor:cocycle}
For $\bP$-a.e.\ $\om$ the corrector satisfies the cocycle property
\[
\chi(s+t,x+y,\om)=\chi(s,x,\om)+\chi(t,y,\tau_{s,x}\om)
\]
for non-negative $s,t\in \mathbb{Q}$.
\end{corollary}

\proof
We have
\begin{align*}
u_\lambda\circ\tau_{s+t,x+y}-u_\lambda=(u_\lambda\circ\tau_{s,x}-u_\lambda)+(u_\lambda\circ\tau_{t,y}-u_\lambda)\circ \tau_{s,x},
\end{align*}
and the claim follows by taking the $L^2(\bP)$-limit along $\lambda'$ on both sides.
\qed

In the following, for any $G:\bZ^d\times\Omega \rightarrow \bR$ we shall write
\[
\| G\|_\om^2:=\sum_{y\sim 0} \mu_{0y}^\om(0) G(y,\om)^2 .
\]
\begin{corollary} \label{cor:Mbracket}
For every $v\in \mathbb{R}^d$ the covariation process of the martingale $M^v:=v\cdot M$ is given by
\begin{align} \label{eq:Mbracket}
\langle v \cdot M\rangle_t=\int_0^t \| v\cdot \Phi \|_{\eta_s}^2 ds,
\intertext{where} \label{def:Phi}
\Phi(x,\om):=x+\chi(0,x,\om).
 \end{align}
\end{corollary}
\proof
First we compute the covariation process of the martingale $M_t^{j,\lambda}$ defined as in the proof of Proposition \ref{prop:constr_corr}. To that aim we define $z_{j,\lambda}(t,x,\om):=x^j+ u_\lambda^{j}(\tau_{t,x}\om)$. Then, by adding \eqref{def_Nk} and \eqref{def_tildeM} we get
\begin{align*}
M_t^{j,\lambda}=z_{j,\lambda}(t,X_t,\om)-z_{j,\lambda}(0,0,\om)-\int_0^t \bar L z_{j,\lambda}(s,X_s,\om) \, ds,
\intertext{where}
\bar Lz_{j,\lambda}(t,x,\om):=D_t z_{j,\lambda}(t,x,\om)+\sum_{y\sim x} \mu_{xy}^\om (t) (z_{j,\lambda}(t,y,\om)-z_{j,\lambda}(t,x,\om)).
\end{align*}
In particular,
\begin{align*}
\bar Lz^2_{j,\lambda}-2 z_{j,\lambda} \bar L z_{j,\lambda}(t,x,\om)&=\sum_{y\in\bZ^d} \mu_{xy}^\om (t) (z_{j,\lambda}(t,y,\om)-z_{j,\lambda}(t,x,\om))^2 \\
&= \sum_{y\in\bZ^d} \mu_{0,y-x}^{\tau_{t,x}\om}(0) \left( y^j-x^j+ [u_\lambda^{j}\circ\tau_{0,y-x}-u_\lambda^{j}]\circ  \tau_{t,x} (\om) \right)^2 \\
&=\sum_{y\in\bZ^d} \mu_{0,y}^{\tau_{t,x}\om}(0) \left( y^j+ [u_\lambda^{j}\circ\tau_{0,y}-u_\lambda^{j}]\circ  \tau_{t,x} (\om) \right)^2,
\end{align*}
so that
\begin{align*}
\langle M_t^{j,\lambda}\rangle_t=\int_0^t \sum_{y\in\bZ^d} \mu_{0,y}^{\eta_s} (0) \left( y^j+ [u_\lambda^{j}\circ\tau_{0,y}-u_\lambda^{j}] (\eta_s) \right)^2 ds,
\end{align*}
and by taking limits on both sides along $\lambda'$, we obtain 
\begin{align*}
\langle M^j\rangle_t=\int_0^t \sum_{y\in\bZ^d} \mu_{0y}^{\eta_s}(0) [y^j+\chi^j(0,y,\eta_s)]^2 ds=\int_0^t \| \Phi^j \|_{\eta_s}^2 \, ds.
 \end{align*} 
For an arbitrary  $v\in \bR^d$ a similar computation gives \eqref{eq:Mbracket}.
\qed

We conclude this section with a convergence result, which will imply the annealed invariance principle. Nevertheless, it will be convenient to complete the proof of Theorem~\ref{annealed-ip} in Section~\ref{sec:pf_main} below.

\begin{proposition} \label{prop:conv_corr1}
We have $t^{-1/2} \chi(t,X_t,\om) \to 0$ in $L^2(\bP^*_{0,0})$ as $t\to \infty$.
\end{proposition}
\proof
Consider an arbitrary fixed $j\in\{1,\ldots,d\}$. Still using the notation in the proof of Proposition \ref{prop:constr_corr} we have for every $t$ and any $\lambda>0$,
\begin{align*}
\chi^j(t,X_t,\om)=M_t^j-X_t^j=M_t^j-M_t^{j,\lambda}+u_\lambda^{j}(\eta_t)-u_\lambda^{j}(\om)-\lambda \int_0^t  u_\lambda^{j}(\eta_s) \, ds,
\end{align*}
and by Cauchy-Schwarz we get
\[
|\chi^j(t,X_t,\om)|^2\leq 3 |M_t^j-M_t^{j,\lambda}|^2+3 |u_\lambda^{j}(\eta_t)-u_\lambda^{j}(\om)|^2 + 3 \lambda^2 \int_0^t  |u_\lambda^{j}(\eta_s)|^2 \, ds.
\]
We argue similarly as in the proof of  Proposition \ref{prop:constr_corr}. Using the fact that $\bP$ is an invariant measure for the environment process $\eta$ we obtain
\begin{align}
\bE E_{0,0}^\om |M_t^j-M_t^{j,\lambda}|^2&\leq 2t \|u_\lambda^{j}-u^{j}\|^2_{\mathcal{H}_1}, \label{est_term1_res} \\
\bE E_{0,0}^\om |u_\lambda^{j}(\eta_t)-u_\lambda^{j}(\om)|^2 &\leq 4 \|u_\lambda^{j}\|^2_{L^2(\bP)}, \\
\bE E_{0,0}^\om  \lambda^2 \int_0^t  |u_\lambda^{j}(\eta_s)|^2 \, ds &\leq t \lambda^2  \|u_\lambda^{j}\|^2_{L^2(\bP)} \label{est_term3_res}.
\end{align}
Choosing $\lambda=t^{-1}$ the claim follows by Proposition \ref{prop:conv_u}.
\qed

\section{Subdiffusive Bound on the Corrector} \label{sec:asymp_corr}

In this section we shall prove the following

\begin{proposition} \label{prop:asymp_corr}
Under the assumptions of Theorem~\ref{main-ip} or Theorem~\ref{main-ip2}, there exists an $\alpha<1/2$ such that
\[
\bE E_{0,0}^\om\left[|\chi(n,X_n,\om)|^2\right]=O(n^{2\alpha}).
\]
\end{proposition}

\subsection{Convergence of the Resolvents}

\begin{proposition} \label{prop:ulambda}
Under the assumptions of Theorem~\ref{main-ip} or Theorem~\ref{main-ip2}, there exists an $\alpha<1/2$ such that for every $j=1,\ldots,d$,
\[
\|u_\lambda^{j}\|_{L^2(\bP)}=O(\lambda^{-\alpha}).
\]
\end{proposition}
Note that while for the annealed FCLT the convergence in Proposition~\ref{prop:conv_corr1} is sufficient, we will need the stronger statement in Proposition~\ref{prop:asymp_corr} for the QFCLT. This difference also appears in the corresponding results on the resolvents $u_\lambda$, (cf.\ Proposition~\ref{prop:conv_u} and Proposition~\ref{prop:ulambda}).
Before we prove Proposition~\ref{prop:ulambda} we will first show how it implies Proposition~\ref{prop:asymp_corr}.

\bigskip {\bf Proof of Proposition \ref{prop:asymp_corr}.}
Similarly to the proof of Proposition~\ref{prop:conv_corr1} we show that for a certain $\lambda$ chosen below depending on $n$ the terms in the right hand side of \eqref{est_term1_res}-\eqref{est_term3_res} are in $O(n^{2\alpha})$.
We shall use similar arguments as in \cite{MW}, in particular cf.\ Lemma~2 and Corollary~4 in \cite{MW}. In a first step we will show that
\begin{align} \label{est_normH1}
\|u_\lambda^{j}-u_{\lambda'}^{j}\|_{\mathcal{H}_1}^2\leq \frac {(\sqrt{\lambda}+\sqrt{\lambda'})^2}{2} \left( \|u_\lambda^{j}\|^2_{L^2(\bP)}+ \|u_{\lambda'}^{j}\|^2_{L^2(\bP)}\right).
\end{align}
Indeed, using the fact that  $u_\lambda^{j}$ solves the resolvent equation \eqref{eq:res} we have
\begin{align*}
\|u_\lambda^{j}-u_{\lambda'}^{j}\|_{\mathcal{H}_1}^2
&=\langle u_\lambda^{j}-u_{\lambda'}^{j}, (-L)(u_\lambda^{j}-u_{\lambda'}^{j}) \rangle_{L^2(\bP)} 
=\langle u_\lambda^{j}-u_{\lambda'}^{j}, -(\lambda u_\lambda^{j}-\lambda' u_{\lambda'}^{j}) \rangle_{L^2(\bP)} \\
&=-\lambda  \|u_\lambda^{j}\|^2_{L^2(\bP)} -\lambda'  \|u_{\lambda'}^{j}\|^2_{L^2(\bP)}+ (\lambda+\lambda') \langle u_{\lambda}^{j},u_{\lambda'}^{j}\rangle_{L^2(\bP)} \\
&\leq 2\sqrt{\lambda \lambda'}  \|u_{\lambda}^{j}\|_{L^2(\bP)}  \|u_{\lambda'}^{j}\|_{L^2(\bP)}+(\lambda+\lambda') \langle u_{\lambda}^{j},u_{\lambda'}^{j}\rangle_{L^2(\bP)} \\
&\leq (\sqrt{\lambda}+\sqrt{\lambda'})^2\|u_{\lambda}^{j}\|_{L^2(\bP)}  \|u_{\lambda'}^{j}\|_{L^2(\bP)},
\end{align*}
which gives \eqref{est_normH1}. In particular, choosing $\lambda_k=2^{-k}$, we get
\begin{align*}
\|u_{\lambda_k}^{j}-u_{\lambda_{k-1}}^{j}\|_{\mathcal{H}_1}^2&\leq \frac {(\sqrt{2^{-k}}+\sqrt{2^{-k+1}})^2}{2} \left( \|u_{\lambda_k}^{j}\|^2_{L^2(\bP)}+ \|u_{\lambda_{k-1}}^{j}\|^2_{L^2(\bP)}\right) \\
&= \frac {(\sqrt{2}+1)^2}{2} \lambda_k  \left( \|u_{\lambda_k}^{j}\|^2_{L^2(\bP)}+ \|u_{\lambda_{k-1}}^{j}\|^2_{L^2(\bP)}\right).
\end{align*}
Let now $k_n$ be the integer $k$ such that $2^{k-1}\leq n <2^k$. Then, we use the elementary estimate $\sqrt{a+b}\leq \sqrt{a}+\sqrt{b}$ for any $a,b\geq 0$ to obtain
\begin{align*}
\|u_{\lambda_{k_n}}^{j}-u^{j}\|_{\mathcal{H}_1}\leq \sum_{m=k_n+1}^\infty \|u_{\lambda_m}^{j}-u_{\lambda_{m-1}}^{j}\|_{\mathcal{H}_1}
\leq c \sum_{m=k_n+1}^\infty  \sqrt{\lambda_m} \left( \|u_{\lambda_m}^{j}\|_{L^2(\bP)}+ \|u_{\lambda_{m-1}}^{j}\|_{L^2(\bP)}\right).
\end{align*}
Recall that $ \|u_{\lambda_m}^{j}\|_{L^2(\bP)}=O(\lambda_m^{-\alpha})$ by Proposition \ref{prop:ulambda}. Therefore, for $n$ large enough
\begin{align*}
\|u_{\lambda_{k_n}}^{j}-u^{j}\|_{\mathcal{H}_1} \leq c  \sum_{m=k_n+1}^\infty \lambda_m^{1/2-\alpha}=c \lambda_{k_n}^{1/2-\alpha}=c n^{\alpha-1/2}.
\end{align*}
Thus, the claim follows by choosing $\lambda_{k_n}$ for $\lambda$ in equation \eqref{est_term1_res}-\eqref{est_term3_res}.
\qed

Recall that $(P_t)_{t\geq 0}$ denotes the transition semigroup of the environment process $\eta$.
\begin{lemma} \label{lem:asymp_Pk}
Under the assumptions of Theorem~\ref{main-ip} or Theorem~\ref{main-ip2}, there exists $\alpha<1/2$ such that for every $j=1,\ldots, d$,
\begin{align*}
\big\| \int_{0}^t P_s V_j \, ds \big\|_{L^2(\bP)} \leq c (1\vee t)^\alpha.
\end{align*}
\end{lemma}
Lemma \ref{lem:asymp_Pk}, which will be proven in the next subsection, immediately implies Proposition~\ref{prop:ulambda}.

 \bigskip {\bf Proof of Proposition \ref{prop:ulambda}.}
Since $u_\lambda^{j}$ is the solution of the resolvent equation \eqref{eq:res},
\begin{align*}
u_\lambda^{j} = \int_0^\infty e^{-\lambda s} P_sV_j \, ds= \lambda \int_0^\infty \int_s^\infty e^{-\lambda t} P_s V_j \, dt \, ds=\lambda \int_0^\infty e^{-\lambda t} \int_0^t P_sV_j \, ds \, dt.
\end{align*}
Hence, by Lemma \ref{lem:asymp_Pk} we get that
$$
\| u_\lambda^{j} \|_{L^2(\bP)} \leq c_1 \lambda \int_0^\infty e^{-\lambda t} (1\vee t)^\alpha \, dt
\leq c_1+ c_1  \lambda \int_1^\infty e^{-\lambda t} t^\alpha \, dt \leq c_1 +c_1 \Gamma(\alpha+1) \lambda^{-\alpha},
$$
which is the claim. 
 \qed

\subsection{A Two-Walk Estimate}
In this subsection we prove Lemma \ref{lem:asymp_Pk}. We shall use techniques from \cite[Section 3]{JR}, \cite[Appendix A]{RS2} and \cite{Mou}. 
Denote by $(X_t)_t$ and $(\tilde X_t)_t$  two independent random walks evolving in the same environment $\om$ both starting from zero. We will write $\bP_{2,x,\tilde x}$ in short for the averaged law of $(X,\tilde X)$ starting in $(x,\tilde x)\in \bZ^d\times \bZ^d$ and $\bE_{2,x,\tilde x}$ for the corresponding expectation, i.e.\ $\bE_{2,x, \tilde x}=\bE\otimes E^\om_{0,x}\otimes E^\om_{0,\tilde x}$. For abbreviation we will write $\bP_{2,x}=\bP_{2,x,0}$ and $\bE_{2,x}=\bE_{2,x,0}$ as well as  $\bP_{2}=\bP_{2,0,0}$ and $\bE_{2}=\bE_{2,0,0}$.
Furthermore, let $(Y_t)_{t\geq 0}$ be the continuous time Markov chain evolving in an environment $\om$ with transition probabilities given by
\[
\pi^\om_{s,x}[Y_t\in A]=P^\om[Y_t\in A | Y_s=x]=\sum_{u,v\in \bZ^d} \indicator_{\{v-u\in A\}} p^\om(s,0;t,u) p^\om(s,x;t,v).
\]
The corresponding expectation will be denoted by $E^{\pi,\om}_{s,x}$.
In particular, note that for every $\om$ the law of $ X_t-\tilde X_t$ induced by $E^\om_x\otimes E^\om_0$ is the same as that of $Y_t$.

\begin{lemma} \label{lem:probY}
For any $0\leq s \leq t$ with $t-s\geq 1$, $y\in \bZ^d$ and any ball $B(x,r)$ we have
\begin{align*}
\pi^\om_{s,y} [Y_t\in B(x,r)]\leq c (t-s)^{-d/2} r^d.
\end{align*}
\end{lemma}
\proof
By the heat kernel estimates in Proposition \ref{prop:hke} we have
\begin{align*}
\pi^\om_{s,y} [Y_t\in B(x,r)]&= \sum_{u,v\in \bZ^d} \indicator_{\{v-u\in B(x,r)\}} p^\om(s,0;t,u) p^\om(s,y;t,v)\\
&=\sum_{u\in \bZ^d} p^\om(s,0;t,u) \sum_{v\in \bZ^d}  \indicator_{\{v\in B(x+u,r)\}} p^\om (s,y;t,v) \\
&\leq c (t-s)^{-d/2} r^d,
\end{align*}
which is the claim. \qed

\subsubsection{Proof of Lemma \ref{lem:asymp_Pk} under Assumptions \ref{ass:time_mixing} and \ref{ass:space_mixing}}
Let $d\geq 3$ and assume that \ref{ass:ergodic}-\ref{ass:space_mixing} hold.
It is enough to show that there exists $\beta>1/2$ such that for every $j=1,\ldots, d$,
\begin{align} \label{eq:asymp_Pk}
 \big\|  P_t V_j \big\|_{L^2(\bP)} \leq c (1\vee t )^{-\beta}.
\end{align}
First note that by definition $V_j(\om)=\mu_{0,e_j}^\om(0)-\mu_{0,-e_j}^\om(0)$, so  by Assumptions \ref{ass:ergodic} and \ref{ass:ellipticity} we have $\bE [ V_j]=0$ and $\|V_j\|_{L^\infty(\bP)}\leq 2C_u$, respectively. In particular, it suffices to prove \eqref{eq:asymp_Pk} for $t\geq 1$. Setting
\[
S_n(f):=\sum_{x\in B(0,n)} f(\tau_{0,x}\om),
\]
we have by the translation invariance of $\bP$
\begin{align*}
\bE[ (S_n(P_tV_j))^2] = \sum_{x,y\in B(0,n)} \bE[P_t V_j(\tau_{0,x-y}\om) P_tV_j(\om)]= \sum_{x,y\in B(0,n)} \bE_{2,x-y}[V_j(\tau_{t,X_t}\om)V_j(\tau_{t,\tilde X_t}\om)].
\end{align*}
Let $\kappa>0$ to be chosen below. Then, for every $z\in \bZ^d$ we use Lemma~\ref{lem:probY} and Assumption~\ref{ass:space_mixing} and obtain
\begin{align} \label{est:X_Xtilde}
 \bE_{2,z}[V_j(\tau_{t,X_t}\om)V_j(\tau_{t,\tilde X_t}\om)]&\leq \bE_{2,z}[V_j(\tau_{t,X_t}\om)V_j(\tau_{t,\tilde X_t}\om)\indicator_{\{|\tilde X_t - X_t|>n^\kappa \}}]+ c \bP_{2,z}[ |Y_t| \leq n^\kappa ] \nonumber \\
 &\leq\bE_{2,z} \left[  \bE_{2,z} \big[V_j(\tau_{t,X_t}\om) V_j(\tau_{t,\tilde X_t}\om) \big| \,X_t, \tilde X_t \big]  \indicator_{\{|\tilde X_t - X_t|>n^\kappa \}} \right] + c t^{-d/2} n^{\kappa d} \nonumber \\
  &=\bE_{2,z} \left[  \bE \big[V_j(\tau_{t,X_t}\om) V_j(\tau_{t,\tilde X_t}\om)]  \indicator_{\{|\tilde X_t - X_t|>n^\kappa \}} \right] + c  t^{-d/2} n^{\kappa d}\nonumber \\
  &\leq c\left( n^{-\kappa p_2} +  t^{-d/2} n^{\kappa d} \right).
\end{align}
Hence,
\begin{align} \label{est_Sn}
n^{-2d}\, \bE[ (S_n(P_tV_j))^2] \leq c\left( n^{-\kappa p_2} +  t^{-d/2} n^{\kappa d}\right).
\end{align}
Next we rewrite the Dirichlet form of the process $\eta$ as
\begin{align} \label{eq:DF_rewritten}
&\langle P_t V_j,P_t V_j \rangle_{\mathcal{H}_1}=\tfrac 1 2 \sum_{y\in\bZ^d} \bE \left[ \mu_{0y}^\om(0) (P_tV_j(\tau_{0,y}\om) -P_tV_j(\om))^2 \right]\nonumber \\
=&\tfrac 1 2 \sum_{y\in\bZ^d} \bE \left[ \mu_{0y}^\om(0) \left(   E^\om_{0,y}[V_j(\tau_{t,X_t}\om)]-E^\om_{0,0}[V_j(\tau_{t,X_t}\om)] \right)\left(   E^\om_{0,y}[V_j(\tau_{t,\tilde X_t}\om)]-E^\om_{0,0}[V_j(\tau_{t,\tilde X_t}\om)] \right) \right]\nonumber \\
=& \tfrac 1 2 \sum_{y\in\bZ^d} \left( \bE_{2,y,y} [ \mu_{0y}^\om(0) V_j(\tau_{t,X_t}\om) V_j(\tau_{t,\tilde X_t}\om) ] -  \bE_{2,y,0} [ \mu_{0y}^\om(0) V_j(\tau_{t,X_t}\om) V_j(\tau_{t,\tilde X_t}\om) ] \right. \nonumber \\
& \left.
- \bE_{2,0,y} [ \mu_{0y}^\om(0) V_j(\tau_{t,X_t}\om) V_j(\tau_{t,\tilde X_t}\om) ] + \bE_{2,0,0} [ \mu_{0y}^\om(0) V_j(\tau_{t,X_t}\om) V_j(\tau_{t,\tilde X_t}\om) ]  \right).
\end{align}
Then, by the time mixing in Assumption \ref{ass:time_mixing} we have
\begin{align*}
\bE_{2,y,y} \left[ \mu_{0y}^\om(0) V_j(\tau_{t,X_t}\om) V_j(\tau_{t,\tilde X_t}\om) \right]&= \bE_{2,y,y}\left[ \bE_{2,y,y} \big[  \mu_{0y}^\om(0) V_j(\tau_{t,X_t}\om) V_j(\tau_{t,\tilde X_t}\om) \big| \,X_t, \tilde X_t \big] \right]  \\
&=\bE_{2,y,y}\left[ \bE \big[  \mu_{0y}^\om(0) V_j(\tau_{t,X_t}\om) V_j(\tau_{t,\tilde X_t}\om) \big] \right] \\
&\leq \bE_{2,y,y} \left[ \bE \big[  \mu_{0y}^\om(0) \big] \bE\big[ V_j(\tau_{t,X_t}\om) V_j(\tau_{t,\tilde X_t}\om) \big] \right] + c t^{-p_1} \\
& =\bE \big[  \mu_{0y}^\om(0) \big] \cdot \bE_{2,0,0}\left[ \bE\big[ V_j(\tau_{t,X_t}\om) V_j(\tau_{t,\tilde X_t}\om) \big] \right]+ c t^{-p_1} \\
&\leq  c \left( n^{-\kappa p_2} +  t^{-d/2} n^{\kappa d} +  t^{-p_1} \right),
\end{align*}
where we also used Assumption \ref{ass:ergodic} in the fourth step and \eqref{est:X_Xtilde} in the last step. The other three terms in \eqref{eq:DF_rewritten}
can be treated similarly, and we obtain that
\begin{align} \label{est:DF}
\langle P_t V_j,P_t V_j \rangle_{\mathcal{H}_1}\leq  c \left( n^{-\kappa p_2} +  t^{-d/2} n^{\kappa d} +  t^{-p_1} \right).
\end{align}
Note that by the ellipticity in Assumption \ref{ass:ellipticity} for any $f\in \mathcal{D}(L)$ the Dirichlet form $\langle f,f \rangle_{\mathcal{H}_1}$ is comparable 
with the Dirichlet form of the environment process associated with a simple random walk on $\bZ^d$. Thus, by Proposition 3.2 in \cite{Mou}, which is a simple consequence of the local Poincar\'{e} inequality on $\bZ^d$, there exists $C_S>0$ such that for any $f\in \mathcal{D}(L)$ and $n\in \mathbb{N}$,
\[
\bE[f(\om)^2]\leq  C_S n^2 \langle f,f \rangle_{\mathcal{H}_1} + \frac 2 {|B(0,n)|^2} \bE[S_n(f)^2].
\]
Combining this with \eqref{est_Sn} and \eqref{est:DF} we get
\begin{align} \label{est_PtVj} 
\bE[P_tV_j(\om)^2]&\leq  C_S n^2 \langle P_tV_j,P_t V_j \rangle_{\mathcal{H}_1} + \frac 2 {|B(0,n)|^2} \bE[S_n(P_t V_j)^2]\nonumber \\
&\leq c  \left( n^{2-\kappa p_2} +  t^{-d/2} n^{2+\kappa d} +  t^{-p_1} n^2 \right).
\end{align}
By Assumption \ref{ass:space_mixing} we have  $p_2> 2d/(d-2)$, so there exists $\delta>0$ such that  $p_2>(1+\delta) 2d/(d-2)$.  Now let 
\begin{align*}
\kappa> \max \left(\frac {2+\frac{4(1+\delta)}{d-2}} {p_2-\frac{2(1+\delta)d}{d-2}} , \frac 1 d \left( \frac{d-2}{p_1-1}-2 \right) \right)
 \intertext{and} 
 \varrho\in \left( \frac 1 {1+\delta} \frac{d/2-1}{\kappa d+2},\frac{d/2-1}{\kappa d+2} \right)  .
\end{align*}
Finally, choosing $n=t^\varrho$ in \eqref{est_PtVj} gives \eqref{eq:asymp_Pk}.
\qed

\subsubsection{Proof of Lemma \ref{lem:asymp_Pk} under Assumptions \ref{ass:time_mixing2} and \ref{ass:space_mixing2}}
Let $d\geq 1$ and assume that \ref{ass:ergodic}-\ref{ass:ellipticity}, \ref{ass:time_mixing2} and \ref{ass:space_mixing2} hold.
Notice first that
\begin{align*} 
\begin{split}
\big\| \int_0^t P_s V_j \, ds \big\|_{L^2(\bP)}^2&=\int_0^t \int_0^t \bE \bigl[ \sum_{x,y\in \bZ^d} V_j(\tau_{r,x}\om)V_j(\tau_{s,y}\om)  p^\om(0,0;r,x) p^\om(0,0;s,y) \bigr] \, dr \, ds\\
&=
\int_0^t \int_0^t \bE \left[ E^\om_{0,0}[V_j(\tau_{r,X_r}\om)] E^\om_{0,0}[V_j(\tau_{s,\tilde X_s}\om)]\right] \, dr \, ds  \\
&=2 \int_0^t \int_0^t \indicator_{\{ r\leq s \}} \bE \left[ E^\om_{0,0}[V_j(\tau_{r,X_r}\om)] E^\om_{0,0}[V_j(\tau_{s,\tilde X_s}\om)]\right] \, dr \, ds 
\end{split}
\end{align*}
and that by definition $V_j(\om)=\mu_{0,e_j}^\om(0)-\mu_{0,-e_j}^\om(0)$, so by Assumption \ref{ass:ergodic} and \ref{ass:ellipticity} we have $\bE [ V_j]=0$ and $\|V_j\|_{L^\infty(\bP)}\leq 2C_u$, respectively. Again it suffices to consider $t\geq 1$.

For any $0\leq r <s \leq t$ with $s-r \geq 1$ we have by Assumption \ref{ass:time_mixing2}
\begin{align*}
\bE E^\om_{0,0}[V_j(\tau_{r,X_r}\om)] E^\om_{0,0}[V_j(\tau_{s,\tilde X_s}\om)]
&=\bE_2 \left[  \bE_2 \big[V_j(\tau_{r,X_r}\om) V_j(\tau_{s,\tilde X_s}\om) \big| \,X_r, \tilde X_s \big]  \right] \\
&=\bE_2 \bE \left[ V_j(\tau_{r,X_r}\om) V_j(\tau_{s,\tilde X_s}\om) \right] \\
&\leq c (s-r)^{-p_1}
\end{align*}
with $p_1>d+1$ if $d\geq 2$ and $p_1>4$ if $d=1$. We fix $\delta\in (1/p_1,1/(d+1))$ if $d\geq 2$ and $\delta\in (1/p_1,1/4)$ if $d=1$ and set $T_t:= t^\delta$ and $L_t:=(1\vee c_1) T_t$ (with constant $c_1$ as in Proposition \ref{prop:hke}). Then,
\begin{align} \label{eq:kl_apart} \begin{split}
\int_0^t \int_0^t   \indicator_{\{ r\leq s \}} \indicator_{\{s-r\geq T_t \}}  \bE\left[ E_{0,0}^\om[V_j(\tau_{r,X_r}\om)] E_{0,0}^\om[V_j(\tau_{s,\tilde X_s}\om)] \right]\, dr \, ds&\leq c t^2 T_t^{-p_1} .
\end{split}
\end{align}
Now we shall consider pairs of times $r$ and $s$ with distance less than $T_t$. We decompose the integral as follows.
\begin{align} \label{eq:kl_close}
&\int_0^t \int_0^t   \indicator_{\{ r\leq s \}}   \indicator_{\{ s-r< T_t \}} \bE \left[ E_{0,0}^\om[V_j(\tau_{r,X_r}\om)] E_{0,0}^\om[V_j(\tau_{s,\tilde X_s}\om)] \right]\, dr \, ds \nonumber \\
\leq &
\int_0^t \int_0^t    \indicator_{\{ r\leq s \}} \indicator_{\{ s-r< T_t \}}  \bE_2 \bigl[V_j(\tau_{r,X_r}\om) V_j(\tau_{s,\tilde X_s}\om) \indicator_{\{|X_r-\tilde X_r| >2 L_t\}} \bigr] \, dr \, ds \nonumber \\
 &+ c T_t\int_0^t  \bP_2[|X_r-\tilde X_r| \leq 2 L_t] \, dr \nonumber \\ 
\leq & \int_0^t \int_0^t   \indicator_{\{ r\leq s \}}  \indicator_{\{ s-r< T_t \}} \bE_2 \bigl[ V_j(\tau_{r,X_r}\om) V_j(\tau_{s,\tilde X_s}\om) \indicator_{\{|X_r-\tilde X_s| >L_t\}}  \indicator_{\{|X_r-\tilde X_r| >2 L_t\}}\bigr] \, dr \, ds\nonumber  \\
&+ c  \int_0^t \int_0^t  \indicator_{\{ r\leq s \}}   \indicator_{\{ s-r< T_t \}} \bP_2[|X_r-\tilde X_s| \leq L_t, \, |X_r-\tilde X_r|> 2 L_t] \, dr \, ds \nonumber  \\
& + c T_t \int_0^t  \bE \pi^\om_{0,0} [|Y_r| \leq 2 L_t] \, dr.
\end{align}
To estimate the first term in \eqref{eq:kl_close} note that conditioned on the event $\{ |X_r-\tilde X_s| >L_t\}$ we have that $V_j(\tau_{r,X_r}\om)$ 
and  $V_j(\tau_{s,\tilde X_s}\om)$ depend only on variables contained in two subsets of $\bZ^d$ with distance $L_t$. Thus, by Assumption \ref{ass:space_mixing2} we obtain
\begin{align} \label{est_term1}
&\bE_2 \left[ V_j(\tau_{r,X_r}\om) V_j(\tau_{s,\tilde X_s}\om) \indicator_{\{|X_r-\tilde X_s| >L_t\}}  \indicator_{\{|X_r-\tilde X_r| >2 L_t\}} \right] \nonumber \\
=& \bE_2 \left[  \bE_2 \big[V_j(\tau_{r,X_r}\om) V_j(\tau_{s,\tilde X_s}\om) \big| \,X_r, \tilde X_r , \tilde X_s \big] \indicator_{\{|X_r-\tilde X_s| >L_t\}}  \indicator_{\{|X_r-\tilde X_r| >2 L_t\}} \right] \nonumber \\
=&\bE_2 \left[  \bE \big[V_j(\tau_{r,X_r}\om) V_j(\tau_{s,\tilde X_s}\om)  \big] \indicator_{\{|X_r-\tilde X_s| >L_t\}}  \indicator_{\{|X_r-\tilde X_r| >2 L_t\}}\right] \nonumber \\
\leq&  c L_t^{-p_2} .
\end{align}

Next we estimate the second term in \eqref{eq:kl_close}.  First we use the Markov property to get
\begin{align*}
 \bP_2[|X_r-\tilde X_s| \leq L_t, \, |Y_r|> 2 L_t]\leq  \bP_2[|\tilde X_s-\tilde X_r | >L_t]=\bP P_{0,0}^\om P^\om_{r,\tilde X_r}[|\tilde X_s-\tilde X_r | >L_t].
\end{align*}
Set $D_i:=\{ y\in \bZ^d:\, 2^i L_t \leq |y-\tilde X_r | \leq 2^{i+1} L_t\}$, $i\geq 0$. Then, noting that $s-r \leq T_t\leq c_1^{-1} L_t$, we use the heat kernel estimates in Proposition  \ref{prop:hke} to obtain
\begin{align*}
P^\om_{r,\tilde X_r}[|\tilde X_s-\tilde X_r | >L_t]&=\sum_{y\in B(\tilde X_r,L_t)^c} p^\om (r,\tilde X_r;s,y)=\sum_{i=0}^\infty \sum_{y\in D_i} p^\om (r,\tilde X_r;s,y) \\
&\leq c \sum_{i=0}^\infty \sum_{y\in D_i}  \exp\big(-c |y-\tilde X_r| \log(|y-\tilde X_r| / (s-r))\big) \\
&\leq c \sum_{i=0}^\infty (2^i L_t)^d   \exp\big(-c 2^i L_t \log(2^iL_t/T_t)\big)\\
 &\leq c \sum_{i=0}^\infty (2^i L_t)^d   \exp\big(-c 2^i L_t \big).
\end{align*}
An elementary computation 
now gives
\begin{align} \label{est_term2}
\bP_2[|X_r-\tilde X_s| \leq L_t, \, |Y_r| > 2 L_t] &\leq c L_t^d  \int_1^\infty \exp(-c L_t u )\, du \leq c  L_t^{d-1} \exp\left(-c L_t \right). 
\end{align}
To estimate the last term in \eqref{eq:kl_close} we use Lemma~\ref{lem:probY} to obtain in the case $d\geq 2$
\begin{align} \label{est_term3}
\int_0^t  \pi^\om_{0,0}[|Y_r| \leq 2 L_t] \, dr \leq c  \int_0^t (1\vee r)^{-d/2} L_t^d \, dr \leq c L_t^d \int_0^t (1\vee r)^{-1} \, dr \leq c \log t \, L_t^d,
\end{align}
and if $d=1$
\begin{align} \label{est_term3d1}
\int_0^t  \pi^\om_{0,0}[|Y_r| \leq 2 L_t] \, dr \leq c \int_0^t (1\vee r)^{-1/2} L_t^d \, dr  \leq c  t^{1/2} \, L_t.
\end{align}

Finally, combining \eqref{eq:kl_apart} and \eqref{eq:kl_close} we get in the case $d\geq 2$ by \eqref{est_term1}, \eqref{est_term2} and \eqref{est_term3}
\begin{align*}
\big\| \int_0^t P_s V_j \big\|_{L^2(\bP)}^2 \, ds \leq c \left(t^2 T_t^{-p_1} +t T_t L_t^{-p_2} + t T_t L_t^{d-1} \exp(-c {L_t})+T_t \log t \, L_t^d \right).
\end{align*}
Analogously, if $d=1$ we obtain by  \eqref{est_term1}, \eqref{est_term2} and \eqref{est_term3d1} that
\begin{align*}
\big\| \int_0^t P_s V_j \big\|_{L^2(\bP)}^2 \, ds \leq c \left(t^2 T_t^{-p_1} +t T_t L_t^{-p_2} + t T_t \exp(-c {L_t})+ t^{1/2}  T_t \, L_t \right).
\end{align*}
The claim follows by our choice of $\delta$, $L_t$ and $T_t$.
\qed

\section{Invariance Principle for $X$} \label{sec:pf_main}
In this section we prove the annealed FCLT in Theorem~\ref{annealed-ip} and the quenched FCLT in Theorem \ref{main-ip} and Theorem \ref{main-ip2}, respectively.  Throughout this section we suppose that Assumptions \ref{ass:ergodic}-\ref{ass:ellipticity} hold. 
The first step to prove a quenched invariance principle for the random walk $X$ is to show that the processes $X^{(\eps)}$ are tight. 
\begin{theorem} \label{t-tight}
Let $T>0$, $r>0$. Then
\begin{align*} 
  & \lim_{R \to \infty} \sup_{0<\eps\leq 1}
    P_{0,0}^\om( \sup_{s \le T} | X^{(\eps)}_s| > R ) \to 0, \\
  & \lim_{\delta \to 0} \limsup_{\eps \to 0}
P_{0,0}^\om( \sup_{|s_1-s_2|  \le \delta, s_i \le T }
 | X^{(\eps)}_{s_2}- X^{(\eps)}_{s_1}| > r )= 0.
\end{align*}
\end{theorem}
\proof  From the heat kernel estimates in Proposition \ref{prop:hke} one can derive tail estimates for the exit times of $X$ from balls (see e.g.\ \cite[Proposition~4.7]{ABDH}).  Then tightness follows by the same arguments as in \cite[Proposition~5.13]{ABDH}.
\qed

\ms For $n \in \bN$ let $\wh X_n=X_n$, and set
\begin{equation}
  \wh X^{(\eps)}_t = \eps \wh X_{\lfloor t/\eps^2 \rfloor}, \q 0< \eps \le 1.
\end{equation}

\begin{lemma}\label{disc-approx}
For any $u>0$,
\begin{equation}
  \lim_{\eps \to 0} P_{0,0}^\om( \sup_{0\le s \le T} 
 |  \wh X^{(\eps)}_s -   X^{(\eps)}_s| > u ) =0.
\end{equation}
\end{lemma}

\proof
This follows from the proof of Theorem~\ref{t-tight} by the same arguments as in \cite[Lemma 4.12]{BD}.
\qed

We will first establish the convergence of the processes  $\wh X^{(\eps)}$; using Lemma~\ref{disc-approx} will then give the convergence of $X^{(\eps)}$ to the same limit. We define
\begin{equation}\label{whM-def}
  \wh M_n =M_n, \qq
 \wh M^{(\eps)}_t =  \eps \wh M_{ \lfloor t/\eps^2 \rfloor}, \q t \ge 0,
\end{equation}
so that 
\begin{equation}\label{xmchi}
  \wh X^{(\eps)}_t =   \eps \wh X_{ \lfloor t/\eps^2 \rfloor}=
\wh M^{(\eps)}_t + \eps \chi(\lfloor t/\eps^2 \rfloor, \eps^{-1} \wh X^{(\eps)}_t ,\om). 
\end{equation}
Thus it is sufficient to prove that the martingale
$\wh M^{(\eps)}$ converges to
a Brownian motion with a certain covariance matrix, and that the second term in \eqref{xmchi} converges to zero
 in $P_{0,0}^\om$-probability for $\bP$-a.a.\ $\om$ (resp.\
 in $\bP^*_{0,0}$-probability) to get the quenched FCLT (resp.\ the annealed FCLT).
For any $G:\bZ^d\times\Omega \rightarrow \bR$ we define
\[
\ol \bE [G]=\sum_{y\sim 0} \bE \left[ \mu^\om_{0y}(0) G(y,\om) \right].
\]

\begin{proposition} \label{P:mconv}
For $\bP$-a.e. $\om$,
the sequence of processes $(\wh M^{(\eps)})$ 
converges in law in the Skorohod topology to a Brownian motion with a non-degenerate
 covariance matrix $\Sigma$ given by $\Sigma_{ij} = \obE \Phi_i \Phi_j$.
\end{proposition}

\proof
We proceed as in \cite{BB}.
Let $v \in \bR^d$ be a unit vector, write as before $\wh M^v_n= v\cdot M_n$, and
let
$$ F_K(\om) = E^0_\om (|\wh M^v_1|^2; |\wh M^v_1| \ge K ). $$
Then $F_K$ is decreasing in $K$, in particular $ \bE [F_K ]\leq \bE [F_0]$.
In the notation of Corollary~\ref{cor:Mbracket} 
$F_0(\om)=\|v \cdot \Phi \|^2_\om$, and so by \eqref{eq:Mbracket}
the covariance process of $\wh M^v$ is
$$ \langle \wh M^v \rangle_n = \int_0^n F_0( \eta_s) \, ds. $$ 
So by the ergodicity of the environment process $\eta$ w.r.t.\ $\bP$ we have $n^{-1} \langle \wh M^v \rangle_n  \to \bE [F_0]$ as $n\to \infty$,
$P_{0,0}^\om$ a.s., for $\bP$-a.a. $\om$.

Using the same arguments as in \cite[Theorem 6.2]{BB} it is
straightforward to check the conditions of the Lindeberg-Feller FCLT
for martingales (see for example \cite[Theorem 3.4.5]{Du}), and deduce
that $v \cdot \wh M^{(\eps)}$ converges to a real-valued
Brownian motion with non-random covariance  $\bE[ \| v\cdot \Phi \|^2_\om ]$,
which can be written as $v\cdot \Sigma v$, where $\Sigma$ is the matrix with coefficients given by
$\Sigma_{ij} = \obE [\Phi_i \Phi_j]$. By the  Cramer-Wold Theorem (see e.g.\ Theorem 3.9.5 in \cite{Du}) 
we get that $\wh M^{(\eps)}$ converges in law to an $\bR^d$-valued Brownian motion with covariance matrix $\Sigma$. 

It remains to show that $\Sigma$ is non-degenerate. By the uniform lower bound on the conductances in Assumption \ref{ass:ellipticity} we have for every unit vector $v\in\bR^d$ that $v\cdot \Sigma v\geq v\cdot \Sigma_{C_l} v$, where $\Sigma_{C_l}$ denotes the non-degenerate covariance matrix of the limiting Brownian motion in the invariance principle for the simple random walk on $\bZ^d$ with constant jump rate $C_l$. Thus, $v\cdot \Sigma v>0$, which implies that $\Sigma$ is positive-definite.
\qed

To conclude the proof of the invariance principles we need to control the corrector function. First we complete the proof of the annealed FCLT.

\bigskip {\bf Proof of Theorem \ref{annealed-ip}.}
Setting $R_n:=\chi(n,X_n,\om)$ we need to show that
\begin{align} \label{conv_Rn_annealed}
n^{-1/2} \max_{k\leq n} |R_k| \rightarrow 0 \qquad \mbox{in $\bP^*_{0,0}$-probability as $n\to \infty$.}
\end{align}
By Proposition~\ref{prop:conv_corr1} we have that  $n^{-1/2}  R_n$ converges to $0$ in $L^2(\bP^*_{0,0})$ and thus in $\bP^*_{0,0}$-probability. By  an elementary property of real convergent sequences, we get \eqref{conv_Rn_annealed}.
\qed

Finally, to complete the proof of the quenched invariance principle we prove
\begin{proposition}
 Let $T>0$. Under the assumptions of Theorem \ref{main-ip} or Theorem \ref{main-ip2}, for $\bP$-a.e.\ $\om$,
\[
\sup_{t\leq T}\eps \chi(\lfloor t/\eps^2 \rfloor, \eps^{-1} \wh X^{(\eps)}_t ,\om)\rightarrow 0 \qquad \mbox{in $P^\om_{0,0}$-probability.}
\]
\end{proposition}

\proof We will proceed as in \cite{RS} applying the theory of ``fractional coboundaries'' of Derriennic and Lin \cite{DL}. Setting $R_n:=\chi(n,X_n,\om)$ we need to show that
\begin{align} \label{conv_Rn}
n^{-1/2} \max_{k\leq n} |R_k| \rightarrow 0 \qquad \mbox{in $P_{0,0}^\om$-probability as $n\to \infty$.}
\end{align}

Let $\tilde \bP$ denote the path measure on $\Omega^{\bN}$ of the random sequence $(\tau_{n,X_n}\om)_{n\geq 0}$ with initial distribution $\bP$, and let $\theta$ be the shift map on the sequence space $\Omega^{\bN}$. By the cocycle property in Corollary~\ref{cor:cocycle} we have $\chi(0,0,\om)=0$ and hence
\begin{align*}
R_n&=\sum_{k=0}^{n-1} \chi(k+1,X_{k+1},\om)-\chi(k,X_k,\om) =\sum_{k=0}^{n-1}h(\tau_{k,X_k}\om, \tau_{k+1,X_{k+1}}\om)
\end{align*}
with $h$ defined as in Remark \ref{rem:defh}. For sequences $\bar \om=(\om^{(i)})_{i\in \bN}$ define $H(\bar \om)=h(\om^{(0)},\om^{(1)})$ and 
\[
\tilde R_n=\sum_{k=0}^{n-1} H\circ \theta^k.
\]
Then $H\in L^2(\tilde\bP)$ and the process $(\tilde R_n)$ has the same distribution under $\tilde \bP$ as the process $(R_n)$ under $\bP\otimes P_{0,0}^\om$. 

By Proposition \ref{prop:asymp_corr} the assumptions of Theorem 2.17 in \cite{DL} are satisfied. We conclude that $H\in (I-\theta)^\gamma L^2(\tilde \bP)$ for any $\gamma\in (0,1-\alpha)$. Since $\alpha<1/2$ there exists such a $\gamma \in (1/2,1-\alpha)$. Then, (i) in Theorem 3.2 in \cite{DL} implies that $n^{-1/2} \tilde R_n$ converges to $0$,  $\tilde\bP$-a.s. Hence, $n^{-1/2}  R_n$ converges to $0$,  $\bP\otimes P_{0,0}^\om$-a.s. In other words, 
 $n^{-1/2}  R_n$ converges to $0$, $P_{0,0}^\om$-a.s., for $\bP$-a.e.\ $\om$, 
which implies \eqref{conv_Rn}.
\qed

\section{Application to Stochastic Interface Models}\label{sec:interface}
In this section we point out a relation between our results and the stochastic dynamic of an interface describing the separation of two pure thermodynamical phases, known  as the Ginzburg Landau $\nabla\phi$ model. We refer to \cite{F} for a survey on these models. The interface is described by a field of height variables $\phi_t(x)$, $x\in \mathbb{Z}^d$, $t\geq 0$, whose stochastic dynamics are given by the following infinite system of stochastic differential equations involving nearest neighbour interaction: 
\begin{align} \label{phi_dyn}
\phi_t(x)=\phi_x-\int_0^t \sum_{y:|x-y|=1} V'(\phi_t(x)-\phi_t(y)) \, dt + \sqrt{2} w_t(x), \qquad x\in \mathbb{Z}^d.
\end{align}
Here $\phi$ is the height of the interface at time $t=0$, $\{w(x), x\in \mathbb{Z}^d\}$ is a collection of independent Brownian motions and the potential $V\in C^2(\mathbb{R},\mathbb{R}_+)$ is even and strictly convex, i.e.\
\begin{align} \label{Vconvex}
c_-\leq V''\leq c_+,
\end{align}
for some $0<c_-<c_+<\infty$. Let for each $r>0$
\[
\sE_r:=\{ \phi\in\bR^{\bZ^d}:\, \sum_x | \phi_x|^2 e^{-2 r |x|} <\infty\}
\]
denote the set of tempered configurations. Then, for every initial value $\phi\in \sE_r$ the SDE \eqref{phi_dyn} admits a unique strong solution $\phi_t \in \sE_r$, $t\geq 0$, see \cite{FS}.  Let $H$ be  the formal Hamiltonian given by
\[
H(\phi)=\tfrac 1 2  \sum_{y:|x-y|=1} V(\phi_x-\phi_y),
\]
then the formal equilibrium measure for the dynamic is given by the Gibbs measure
\[
\frac 1 Z \exp(-H(\phi)) \prod_x d\phi_x. 
\]
This can be made rigorous for the corresponding dynamic on a finite box. In dimension $d\geq 3$  Gibbs measures for the $\phi$-field on the whole lattice can be constructed by taking the thermodynamical  limit, cf.\ Section 4.5 in \cite{F}. More precisely, for every $h\in \bR$ there exists a shift-invariant and ergodic $\phi$-Gibbs measure $m_h$ with mean $h$, i.e.\
\[
\int \phi_x \, m_h(d\phi)=h, \qquad x\in\bZ^d.
\]
These measures are also reversible and ergodic for the SDE \eqref{phi_dyn}. We denote by  $\bP_{m_h}$ the law of the process $\phi_t$ started under the equilibrium distribution $m_h$ (and by $\bE_{m_h}$ the corresponding expectation).

Next we consider discrete gradients, i.e.\  height differences of the form $\nabla_b \phi=\phi_{y_b}-\phi_{x_b}$ for any bond $b=\{x_b,y_b\} \in E_d$. 
Then, as a vector field $\nabla \phi$ has zero curl in the sense that
\[
\sum_{b\in\sC}\nabla_b\phi=0
\]
for every closed loop $\sC$, i.e.\  the bonds $\{x_i,x_{i+1}\}$ of a sequence of $x_0,\ldots,x_n$ in $\bZ^d$ satisfying $x_0=x_n$ and $|x_i-x_{i-1}|=1$ for $i\in\{1,\ldots n\}$. Let $\sX$ be the subset of $\bR^{E_d}$, whose elements have zero curl, and let for $r>0$
\[
\sX_r=\{ \eta\in \bR^{E_d}: \, \eta_b=\nabla_b \phi \mbox{ for some } \phi\in\sE_r \}
\]
be the subset of tempered gradients. Note that the drift term in the SDE \eqref{phi_dyn} can be rewritten as
\[
-\sum_{y:|x-y|=1} V'(\phi(x)-\phi(y))=\sum_{b:\, x_b=x} V'(\nabla_b \phi).
\]
Then, for each initial $\nabla\phi \in \sX_r$, the gradient process $(\nabla_b\phi_t, b\in E_d, t\geq 0)$ is the unique strong solution of the SDE
\begin{align*}
\nabla_b\phi_t=\nabla_b\phi-\int_0^t \left( \sum_{b': x_{b'}=x_b} V'(\nabla_{b'}\phi_s)-\sum_{b': x_{b'}=y_b} V'(\nabla_{b'}\phi_s)   \right) \, ds+\sqrt 2 \nabla_b w_t, \qquad b\in E_d,  
\end{align*}
where $ \nabla_b w_t=w_t(y_b)-w_t(x_b)$, see again \cite{FS}. Also it has been shown  in \cite[Theorems 3.1 and 3.2]{FS} that in any lattice dimension $d\geq 1$, given any $u\in \bR^d$, there exists a unique shift invariant ergodic $\nabla \phi$-Gibbs measure $\tilde m_u$ on $\sX_r$ satisfying
\[
\int_{\sX} \eta_{0,e_i} \, d\tilde m_u=u_i, 
\]
for every $i=1,\ldots, d$.  
Here $u$ is the tilt and $\tilde m_u$ the $u$-tilted measure.  Moreover, $\tilde m_u$ is known to be an invariant reversible and ergodic measure for the gradient process $\nabla\phi_t$ (\cite[Proposition 3.1]{FS}). 

Our aim is to investigate the decay of the space-time correlation functions  under the equilibrium Gibbs measures. The idea -- originally from Helffer and Sj\"ostrand \cite{HS} -- is to describe the correlation functions in terms of a certain random walk in dynamic random environment (cf.\ also \cite{DD, GOS, BG}).
Let $(X_t)_{t\geq 0}$ be the random walk on $\bZ^d$ with jump rates given by the random dynamic conductances
\begin{align*} 
\mu^{ \nabla \phi}_{b}(t):=V''(\nabla_b \phi(t))=V''(\phi_{y_b}(t)-\phi_{x_b}(t)), \qquad b=\{x_b,y_b\}\in E_d.
\end{align*}
Since $V$ is even, the jump rates are symmetric, i.e.\ $\mu^{\nabla\phi}_{x_b,y_b}(t)=\mu^{\nabla\phi}_{y_b,x_b}(t)$. Further, let $p^{\nabla\phi}(s,x;t,y)$, $x,y\in \bZ^d$, $s\leq t$, denote the transition densities of the random walk $X$.

\begin{theorem}
\begin{enumerate}
\item[i)] Let $d\geq 3$ and let $m_h$ be any ergodic $\phi$-Gibbs measure. Then, the environment $\mu^{\nabla\phi}$ started under $m_h$ satisfies Assumptions \ref{ass:ergodic}-\ref{ass:ellipticity}. Moreover,  $\mu^{\nabla\phi}$ also satisfies Assumptions \ref{ass:time_mixing} and \ref{ass:space_mixing}  if $d\geq 6$. 

\item[ii)] Let $d\geq 1$ and let $\tilde m_u$ be any ergodic $\nabla\phi$-Gibbs measure. Then, the environment $\mu^{\nabla\phi}$ started under $\tilde m_u$ satisfies Assumptions \ref{ass:ergodic}-\ref{ass:ellipticity}. Moreover,  $\mu^{\nabla\phi}$ also satisfies Assumptions \ref{ass:time_mixing} and \ref{ass:space_mixing}  if $d\geq 5$. 
\end{enumerate}
\end{theorem}

\proof
Assumption \ref{ass:ergodic} is immediate from the ergodicity of the Gibbs measures $m_h$ and $\tilde m_u$, respectively. Assumption \ref{ass:stoch_cont} is clear from the pathwise continuity of $\phi_t$ and $\nabla \phi_t$ and the strict convexity of $V$ in \eqref{Vconvex}  guarantees the ellipticity in Assumption \ref{ass:ellipticity}. 

By Theorem 6.1 in \cite{DD}  the time-covariance under the $\phi$-Gibbs measure $m_h$ decays  polynomially with order $d/2-1$ and the space-covariance decays polynomially with order $d-2$. Hence,  Assumptions \ref{ass:time_mixing} and \ref{ass:space_mixing} hold if $d/2-1>1$ and $d-2>2d/(d-2)$ which is the case for $d\geq 6$.

On the other hand, by Theorem 6.2 in \cite{DD}  the time-covariance for $\nabla \phi$ has polynomial decay of order $d/2$ and the space-covariance has polynomial decay of order $d$ . We have $d/2>1$ and $d>2d/(d-2)$ if $d\geq 5$.
\qed

We combine now the  Helffer-Sj\"ostrand representation and the local limit theorem in Theorem \ref{thm:llt} to get a scaling limit for the space-time covariation of the $\phi$-field.

\begin{theorem}
 Let $d\geq 3$ and let $m_h$ be any ergodic $\phi$-Gibbs measure. Then, there exist a non-degenerate covariance matrix $\Sigma$ such that
 \begin{align*}
N^{d+2} \cov_{m_h}(\phi_0(0), \phi_{N^2t}(\lfloor Ny\rfloor)) \rightarrow \int_0^\infty k^{(\Sigma)}_{t+s}(y)\, ds, \qquad \text{as $N\to \infty$,}
\end{align*}
where $k_t$ is the Gaussian kernel with diffusion matrix $\Sigma$ in \eqref{def_k}.
\end{theorem}

\proof By the Helffer-Sj\"ostrand representation (cf.\ equation (6.10) in \cite{DD}) we have
\[
\cov_{m_h}(\phi_0(0), \phi_t(y))=\int_0^\infty \mathbb{E}_{m_h} \left[ p^{\nabla\phi}(0,0;t+s,y) \right] ds.
\]
Using the annealed local limit theorem in Theorem \ref{thm:llt} i) we obtain 
\begin{align*}
N^{d+2} \cov_{m_h}(\phi_0(0), \phi_{N^2t}\left(\lfloor Ny \rfloor)\right)=&N^d \int_0^\infty \mathbb{E}_{m_h}\left[ p^{\nabla\phi}\big(0,0;N^2(t+s),\lfloor Ny\rfloor\big) \right] ds \\
\rightarrow & \int_0^\infty k^{(\Sigma)}_{t+s}(y)\, ds
\end{align*}
as $N\to \infty$, which is the claim.
\qed

Ultimately, we would like to derive an analogous scaling limit for the space-time covariance of the gradient process $\nabla \phi_t$, see also the discussion in \cite[Section 6]{BG}. However, what is still missing until now is a local limit theorem for the gradient of the heat kernel.

%

\sm SA: Institut f\"ur Angewandte Mathematik \\
Rheinische Friedrich-Wilhelms-Universit\"at Bonn \\
Endenicher Allee 60, 53115 Bonn, Germany. \\
andres@iam.uni-bonn.de

\end{document}